\documentclass[preprint,12pt,numbers,sort&compress,nopreprintline]{elsarticle}
\usepackage{amsmath,amssymb,amsthm,mathtools,microtype}
\usepackage{geometry}
\usepackage[hidelinks]{hyperref}
\numberwithin{equation}{section}
\newtheorem{theorem}{Theorem}[section]
\newtheorem{lemma}[theorem]{Lemma}
\newtheorem{proposition}[theorem]{Proposition}
\newtheorem{corollary}[theorem]{Corollary}
\theoremstyle{definition}

\theoremstyle{remark}
\newtheorem{remark}[theorem]{Remark}
\newcommand{\R}{\mathbb R}
\newcommand{\Sph}{\mathbb S}
\newcommand{\GG}{G}
\newcommand{\HH}{\mathcal H}
\newcommand{\Ph}{\Phi_m}
\newcommand{\vol}{\operatorname{vol}}
\newcommand{\one}{\mathbf 1}
\newcommand{\dd}{\mathrm d}
\newcommand{\doi}[1]{\href{https://doi.org/#1}{\nolinkurl{doi:#1}}.}
\journal{Advances in Mathematics}
\begin{document}
\begin{frontmatter}
\title{Fixed and periodic points of bi-dual intersection operators}
 \author[aff1]{Youjiang Lin}
 \affiliation[aff1]{organization={School of Mathematical Sciences, Hebei Normal University},city={Shijiazhuang},country={China}}
\begin{abstract}
We classify the fixed and periodic points, up to positive dilation, of the
operators $\mathbb I_m=\mathbb I_{m,m+1}$ in the indexed bi-dual intersection
family associated with centro-section integrals. For $n\ge3$ and
$1\le m\le n-2$, a star body satisfies $\mathbb I_m^{\,r}K=cK$ for some
$c>0$ and an integer $r\ge1$ if and only if it is an origin-centered
Euclidean ball. No convexity, origin symmetry, or initial differentiability
is assumed. An exact Grassmannian graph representation gives a decreasing
weight and yields a strictly positive linear increment under continuous
Steiner symmetrization unless the body is already symmetric. At a fixed
point, a uniform radial Taylor estimate makes the same volume-preserving
increment quadratic. Sobolev estimates for the spherical Radon transform
and its adjoint supply the required regularity. Independently, a strict
Lyapunov inequality reduces every periodic point up to dilation to a fixed
point. We also prove automatic positivity and smoothness for bounded
nonnegative solutions and local geometric convergence of normalized
iterates near the ball.
\end{abstract}
\end{frontmatter}
\section{Introduction}\label{sec:introduction}

Moments of central section volumes provide a common setting for integral
geometry, affine inequalities, and Minkowski problems. Their first variations
also produce nonlinear operators on star bodies, and hence a further question:
which bodies retain their shape under these operators, or return to the same
shape after finitely many iterations? The identification of extremizers in a
sharp inequality does not, by itself, classify all stationary points or periodic
orbits. In this paper we address this distinction for a member of the bi-dual
intersection family introduced by Cai, Leng, Wu and Xi. We classify its fixed
and periodic points among all star bodies, without convexity, symmetry, or
smoothness assumptions beyond a positive continuous radial function.

The relevant variational framework was developed in two works of
Cai--Leng--Wu--Xi. In \cite{CLWX2025}, they obtained affine dual curvature
measures by differentiating dual affine quermassintegrals along logarithmic
families of Wulff shapes. They established existence results for the associated
even Minkowski problem, introduced an affine-covariant bi-dual intersection
operator, and proved an associated affine isoperimetric inequality. Their
subsequent paper \cite{CLWX2026} considered the centro-section integrals
\begin{equation}\label{intro:centro-moment}
 \widetilde\Psi_{m,q}(K)
 =\int_{G(n,m)}|K\cap\xi|_m^{\,q}\,d\nu_{n,m}(\xi),\qquad q\ne0,
\end{equation}
where $\nu_{n,m}$ is Haar probability on the Grassmannian of $m$-dimensional
linear subspaces of $\mathbb R^n$. They constructed the corresponding
centro-section measures, proved their variation formulas, and established
existence under subspace concentration conditions, with refined conditions
in the hyperplane case. The resulting indexed bi-dual operators are the
starting point of the present work. Our problem concerns their stationary
shapes and iteration, rather than the prescription of an arbitrary spherical
measure.

\subsection{The indexed operator and the main theorem}

Let $\mathcal S_o^n$ denote the star bodies in $\mathbb R^n$: a set
$K\in\mathcal S_o^n$ has the representation
\[
 K=\{tu:u\in\mathbb S^{n-1},\ 0\le t\le\rho_K(u)\},
 \qquad \rho_K\in C(\mathbb S^{n-1}),\quad \rho_K>0.
\]
Write $\sigma_j=|\mathbb S^{j-1}|$ and $\omega_j=|B_2^j|=\sigma_j/j$.
We use the same Radon-transform normalization as
\cite[Section 2.3]{CLWX2025,CLWX2026}:
\begin{align}
 (\mathcal R_m f)(\xi)
 &=\int_{\mathbb S^{n-1}\cap\xi} f(v)\,d\mathcal H^{m-1}(v),\notag\\
 (\mathcal R_m^* g)(u)
 &=\frac{\sigma_m}{\sigma_n}
   \int_{\xi\ni u}g(\xi)\,d\nu_{u,m}(\xi).
 \label{intro:radon}
\end{align}
Here $\nu_{u,m}$ is the invariant probability on the $m$-planes containing $u$.
Thus $\mathcal R_m^*$ is the adjoint for Grassmann probability and
unnormalized spherical area, and
$|K\cap\xi|_m=\mathcal R_m(\rho_K^m/m)(\xi)$.

The indexed bi-dual $m$-intersection body of
\cite[Eq.~(1.4)]{CLWX2026} is given by
\begin{equation}\label{intro:indexed-operator}
 \rho_{\mathbb I_{m,q}K}^{\,n-m}
 =\mathcal R_m^*\bigl[|K\cap\cdot|_m^{\,q-1}\bigr],\qquad 1\le m<n,\quad q\in\mathbb R.
\end{equation}
The operator originally introduced in \cite[Eq.~(1.13)]{CLWX2025} is the
member $q=n$; its $SL(n)$ covariance is proved in
\cite[Proposition 6.1]{CLWX2025}. We study instead the specialization
$q=m+1$. Throughout the paper, we abbreviate
\begin{equation}\label{intro:operator}
 \mathbb I_m:=\mathbb I_{m,m+1},\qquad
 \rho_{\mathbb I_mK}^{\,n-m}
 =\mathcal R_m^*\!\left[
       \left(\mathcal R_m(\rho_K^m/m)\right)^m\right].
\end{equation}
In particular, for $m\le n-2$, our $\mathbb I_m$ is not the affine operator
denoted by $I_m$ in \cite{CLWX2025}. The indices $m+1$ and $n$ coincide
only when $m=n-1$.

The right-hand side of \eqref{intro:operator} is positive, continuous, and
even. Consequently, $\mathbb I_m$ maps star bodies to origin-symmetric
star bodies. It satisfies
\[
 \mathbb I_m(\lambda K)
 =\lambda^{m^2/(n-m)}\mathbb I_mK,\qquad \lambda>0.
\]
A fixed point or periodic point in this paper is understood \emph{up to
positive dilation}. A superscript on $\mathbb I_m$ denotes composition.

\begin{theorem}[Fixed and periodic points]\label{intro:main}
Let $n\ge3$ and $1\le m\le n-2$. For every $K\in\mathcal S_o^n$ and every
integer $r\ge1$,
\begin{equation}\label{intro:classification}
 \mathbb I_m^{\,r}K=cK\ \text{for some }c>0
 \quad\Longleftrightarrow\quad
 K=RB_2^n\ \text{for some }R>0.
\end{equation}
\end{theorem}

The theorem is global in both the input body and the period.
No closeness to a ball is required, and neither convexity nor initial origin
symmetry is assumed. The line-section case $m=1$ has a direct elementary
proof; the principal geometric and regularity arguments concern
$2\le m\le n-2$.

The endpoint is genuinely different. For the classical intersection operator
$I$, defined by $\rho_{IK}(u)=|K\cap u^\perp|_{n-1}$, our normalization gives
\begin{equation}\label{intro:endpoint}
 \mathbb I_{n-1}K=\mathbb I_{n-1,n}K
 =\frac{n-1}{\sigma_n}I^2K.
\end{equation}
Milman, Shabelman and Yehudayoff \cite[Theorem 1.1]{MSY2025} proved that
$I^2K=cK$, for $n\ge3$, characterizes origin-centered ellipsoids.
Thus the fixed-point conclusion of Theorem~\ref{intro:main} cannot be
extended to $m=n-1$ with balls in place of ellipsoids.

\subsection{Variational meaning and relation to earlier rigidity results}

Set
\begin{equation}\label{intro:phi}
 \Phi_m(K):=\frac1{m+1}\widetilde\Psi_{m,m+1}(K)
 =\frac1{m+1}\int_{G(n,m)}|K\cap\xi|_m^{m+1}\,d\nu_{n,m}(\xi).
\end{equation}
Regarded as a functional of the radial function, its derivative in a
continuous direction $h$ is
\begin{equation}\label{intro:variation}
 D\Phi_m(\rho_K)[h]
 =\int_{\mathbb S^{n-1}}
   \rho_{\mathbb I_mK}^{\,n-m}\rho_K^{m-1}h\,d\mathcal H^{n-1}.
\end{equation}
Since $DV(\rho_K)[h]=\int\rho_K^{n-1}h\,d\mathcal H^{n-1}$, the equation
$\mathbb I_mK=cK$ is exactly the positive radial Euler--Lagrange equation
for $\Phi_m-c^{n-m}V$. Equivalently, it describes the stationary points
of the scale-invariant quotient
\[
 \frac{\widetilde\Psi_{m,m+1}(K)^{1/(m+1)}}{V(K)^{m/n}}.
\]
The homogeneities determine the multiplier: at any such stationary point,
$m(m+1)\Phi_m(K)=n c^{n-m}V(K)$.

This interpretation fits the centro-section measures of
\cite[Definition 3.1]{CLWX2026}. For a convex body $K$ containing the origin in its interior,
those measures are
\[
 \widetilde A_{m,q}(K,\eta)
 =\int_{\alpha_K^*(\eta)}
   \rho_K(u)^m\rho_{\mathbb I_{m,q}K}(u)^{n-m}\,d\mathcal H^{n-1}(u),
\]
where $\alpha_K^*$ is the reverse radial Gauss image. In particular, a
convex fixed point of \eqref{intro:operator} satisfies
\begin{equation}\label{intro:cone-relation}
 \widetilde A_{m,m+1}(K,\cdot)
 =n c^{n-m}V_K(\cdot),
\end{equation}
where $V_K(\eta)=n^{-1}\int_{\alpha_K^*(\eta)}\rho_K^n\,d\mathcal H^{n-1}$
is its cone-volume measure. This is an implication of the pointwise fixed-point equation; no converse
for arbitrary nonsmooth convex bodies is needed here. The problem here nevertheless differs from the prescribed-measure
existence problems of \cite{CLWX2025,CLWX2026}: our equation is a radial
stationarity condition, and the body is unknown on both sides of
\eqref{intro:cone-relation}.

The distinction between stationarity and extremality is essential. The
Busemann--Straus--Grinberg inequality, in its extension to star bodies
\cite{Grinberg1991,Gardner2007}, and monotonicity of probability-space norms give
\begin{equation}\label{intro:moment-bound}
 \widetilde\Psi_{m,m+1}(K)^{1/(m+1)}
 \le \widetilde\Psi_{m,n}(K)^{1/n}
 \le \omega_m\left(\frac{V(K)}{\omega_n}\right)^{m/n}.
\end{equation}
Thus a ball maximizes the quotient above. Such a bound alone does not rule
out other positive stationary configurations. The fixed-point part of
Theorem~\ref{intro:main} excludes precisely those additional configurations;
it is not obtained by assuming equality in \eqref{intro:moment-bound}.

For the classical intersection body operator, introduced by
Lutwak \cite{Lutwak1988}, Fish, Nazarov, Ryabogin and Zvavitch
\cite{FNRZ2011} proved local attraction to the ball in the Banach--Mazur
metric for origin-symmetric star bodies. The global classification of
$I^2K=cK$ and $IK=cK$ was established in \cite{MSY2025}, using a variational
argument and continuous Steiner symmetrization of Lipschitz star bodies.
Lin and Xiong \cite[Theorem 1.1]{LinXiong2025} proved the corresponding
ball classification for $I_i^2K=cK$, $1\le i\le n-2$, where
\begin{equation}\label{intro:lower-order}
 \rho_{I_iK}(u)=\frac1{n-1}
 \int_{\mathbb S^{n-1}\cap u^\perp}\rho_K(v)^i\,d\mathcal H^{n-2}(v).
\end{equation}
These are different lower-order constructions: in
\eqref{intro:lower-order} the hyperplane transform is fixed and the radial
power varies, whereas \eqref{intro:operator} uses the $m$-plane transform,
its Grassmannian adjoint, and two nonlinear powers. Their result therefore
does not directly classify the operators considered here. Both arguments
build on the continuous-symmetrization framework of \cite{MSY2025}; the
additional mechanism in our setting is a strictly decreasing weight in
Grassmannian graph coordinates.

\subsection{A strict Lyapunov functional}

The exclusion of periodic orbits is independent of the geometric fixed-point
classification. Let $\sigma=\mathcal H^{n-1}/\sigma_n$ be sphere probability,
and put
\[
 \mathcal A_m=\frac1{\sigma_m}\mathcal R_m,
 \qquad
 \mathcal A_m^*=\frac{\sigma_n}{\sigma_m}\mathcal R_m^*.
\]
These are adjoint probability averages. Define
\begin{equation}\label{intro:energy}
 \mathcal E_m(K)
 =\frac{\|\mathcal A_m(\rho_K^m)\|_{L^{m+1}(G(n,m),\nu_{n,m})}}
 {\|\rho_K^m\|_{L^{n/m}(\mathbb S^{n-1},\sigma)}}.
\end{equation}
It is invariant under dilation and differs from the quotient in
\eqref{intro:moment-bound} only by a dimensional constant.

\begin{theorem}[Strict monotonicity under iteration]\label{intro:lyapunov}
Let $n\ge2$ and $1\le m<n$. For every $K\in\mathcal S_o^n$,
\[
 \mathcal E_m(\mathbb I_mK)\ge\mathcal E_m(K),
\]
with equality if and only if $\mathbb I_mK=cK$ for some $c>0$.
Consequently, every periodic point up to dilation is already a fixed point
up to dilation.
\end{theorem}

The proof is a nonlinear power-step argument, related to Boyd's method
\cite{Boyd1974}. With $f=\rho_K^m$ and $p=n/m$, the next radial function
raised to the power $m$ is proportional to
\[
 f_+=\bigl\{\mathcal A_m^*[(\mathcal A_mf)^m]\bigr\}^{1/(p-1)}.
\]
Two applications of H\"older's inequality give the claimed increase.
Equality in the second one forces $f_+$ to be proportional to $f$.
Since a dilation-invariant monotone functional has equal initial and terminal
values on a projective periodic orbit, equality already holds at the first
step. No compactness of the orbit or prior classification of its elements
is needed.

\subsection{The proof mechanism}

The fixed-point proof has an analytic preparation and a geometric core.
First, for $m\ge2$, the equation itself implies smoothness. With
$N_m=\mathcal A_m^*\mathcal A_m$, the normal operator has a positive
integral kernel and acts on spherical harmonics of degree $2k$ by
\[
 \lambda_{2k}
 =\frac{(1/2)_k((n-m)/2)_k}{(m/2)_k((n-1)/2)_k}
 \asymp(1+k)^{-(m-1)}.
\]
This gives a gain of $(m-1)/2$ Sobolev orders for each transform.
The estimate for the adjoint is established on its full Grassmannian
domain, which is essential after the nonlinear power has been applied.
Together with independent supremum bounds and Sobolev product estimates,
these gains can be iterated through the fixed-point equation. In fact,
Theorem~\ref{thm:smooth} proves positivity and smoothness for every nonzero
nonnegative bounded solution of the corresponding equation.

For the geometric argument we use the continuous Steiner construction
for Lipschitz star bodies developed in \cite[Sections 4--5]{MSY2025}.
Fix $u\in\mathbb S^{n-1}$ and represent almost every $m$-plane as
\[
 \xi(E,a)=\{y+\langle a,y\rangle u:y\in E\},
 \qquad E\in G(u^\perp,m),\quad a\in E.
\]
Writing $Q_K(E,a)=\int_E\mathbf1_K(y+\langle a,y\rangle u)\,dy$, we prove
\begin{equation}\label{intro:graph}
 \Phi_m(K)=\frac{A_{n,m}}{m+1}
 \int_{G(u^\perp,m)}\int_E
 \frac{Q_K(E,a)^{m+1}}{(1+|a|^2)^{(n-m-1)/2}}\,da\,d\nu(E),
\end{equation}
where $A_{n,m}=\Gamma(n/2)/(\pi^{m/2}\Gamma((n-m)/2))$.
The weight results from the section Jacobian and the Haar density;
it is strictly decreasing in $|a|$ for $m\le n-2$ and becomes constant
at $m=n-1$.

Expanding the integer power in \eqref{intro:graph} produces $m+1$
one-dimensional fiber constraints on an $m$-dimensional slope variable.
A Gaussian-mixture argument and Pr\'ekopa's marginal theorem give
monotonicity along the finite-interval Steiner flow. The further step
needed for rigidity is quantitative at first order: if a Lipschitz star
body is not Steiner symmetric in an admissible direction $u$, then
\begin{equation}\label{intro:strict-growth}
 \Phi_m(S_u^tK)-\Phi_m(K)\ge\kappa_u t
 \qquad(0<t\le\varepsilon_u),\qquad\kappa_u>0.
\end{equation}
We establish strictness at a configuration with a noncentered active
constraint and then retain it uniformly on an open set of parameters.
The last constraint remains redundant near a central fiber. This
positive-measure localization, rather than a pointwise strictness
assertion, yields \eqref{intro:strict-growth}.

At a fixed point, the first variation instead cancels that of volume.
A uniform second-order estimate and the radial displacement bound for
the same volume-preserving flow imply
\begin{equation}\label{intro:quadratic}
 \bigl|\Phi_m(S_u^tK)-\Phi_m(K)\bigr|\le C_ut^2.
\end{equation}
No derivative of the Steiner path with respect to time is required.
Equations \eqref{intro:strict-growth} and \eqref{intro:quadratic} are
incompatible for small positive $t$. Hence the body is reflection
symmetric in almost every admissible direction, and continuity extends
this symmetry to all directions. It is therefore a ball.

The same spectrum gives a complementary local result. On the mean-zero
subspace of $H^s(\mathbb S^{n-1})$, $s>(n-1)/2$, the derivative of the
mean-normalized map at the ball has norm $m/(n-1)<1$. Nearby normalized
iterates converge geometrically in $H^s$; see
Proposition~\ref{prop:local-gap}. Neither this local result nor the
periodic-point classification asserts convergence of every normalized orbit.

Section~\ref{sec:preliminaries} fixes the measure conventions and proves
the radial variational estimates. Section~\ref{sec:regularity} establishes
automatic regularity. The graph and fiber formulas are proved in
Section~\ref{sec:graph}, and the weighted monotonicity argument is developed
in Section~\ref{sec:weighted-monotonicity}. Section~\ref{sec:weighted-strict}
proves strict growth and the fixed-point classification.
Section~\ref{sec:periodic} establishes the Lyapunov inequality and completes
the proof of Theorem~\ref{intro:main}. Section~\ref{sec:local} treats local
dynamics. Two appendices record an independent singular-integral
representation and an elementary maximum--minimum argument.

\section{Preliminaries and radial stationarity}\label{sec:preliminaries}

Throughout, $G(n,m)$ carries Haar probability $\nu_{n,m}$, while
$du=d\mathcal H^{n-1}(u)$ denotes unnormalized spherical area and
$d\sigma=du/\sigma_n$ denotes spherical probability. On a section sphere,
Hausdorff measure has its indicated dimension. We use the transforms
\eqref{intro:radon} and the abbreviation
$\mathbb I_m=\mathbb I_{m,m+1}$ in \eqref{intro:operator}.
Unless a narrower range is specified, $n\ge2$ and $1\le m<n$.

\subsection{Incidence averages and basic mapping properties}

The incidence space
\[
 \mathcal F_{n,m}=\{(u,\xi)\in\Sph^{n-1}\times G(n,m):u\in\xi\}
\]
has a unique orthogonally invariant probability measure. Its marginals
are $\sigma$ and $\nu_{n,m}$. Its conditional measures are the uniform
probability on $\Sph^{n-1}\cap\xi$ and the invariant probability
$\nu_{u,m}$ on the planes containing $u$. Thus the probability averages
\begin{equation}\label{eq:normalized-radon}
 \begin{split}
 (\mathcal A_m f)(\xi)
 &=\frac1{\sigma_m}\int_{\Sph^{n-1}\cap\xi}f(v)\,d\mathcal H^{m-1}(v),\\
 (\mathcal A_m^*g)(u)
 &=\int_{\xi\ni u}g(\xi)\,d\nu_{u,m}(\xi)
 \end{split}
\end{equation}
satisfy
\begin{equation}\label{eq:normalized-adjoint}
 \int_{G(n,m)}(\mathcal A_m f)g\,d\nu_{n,m}
 =\int_{\Sph^{n-1}}f\mathcal A_m^*g\,d\sigma.
\end{equation}
Tonelli proves this identity for nonnegative inputs; integrability then
gives the signed version. For measurable inputs the averages are defined
almost everywhere. A change on a null set changes neither output class,
as follows by applying the same incidence identity to its indicator.

The two normalizations are related by
\begin{equation}\label{eq:normalization}
 \mathcal R_m=\sigma_m\mathcal A_m,
 \qquad \mathcal R_m^*=\frac{\sigma_m}{\sigma_n}\mathcal A_m^*.
\end{equation}
In particular,
\begin{equation}\label{eq:operator-normalized}
 \rho_{\mathbb I_mK}^{n-m}
 =C_{n,m}\mathcal A_m^*[(\mathcal A_m(\rho_K^m))^m],
 \qquad C_{n,m}=\frac{\sigma_m}{\sigma_n}\omega_m^m.
\end{equation}
The distinction between \eqref{eq:normalized-adjoint} and the adjoint
identity for unnormalized area will be maintained in every variation
formula.

\begin{lemma}[Basic mapping properties]\label{lem:averaging}
Both $\mathcal A_m$ and $\mathcal A_m^*$ preserve constants and
positivity, and are contractions on $L^2$ and on $L^\infty$, with their
respective probability measures. They map continuous functions to
continuous functions. In particular, $\mathbb I_m$ maps star bodies to
origin-symmetric star bodies.
\end{lemma}
\begin{proof}
The $L^\infty$ statements follow from averaging. Jensen's inequality
followed by incidence Fubini gives the $L^2$ statements. For a rotation
$U$, equivariance and the supremum norm bound give, for example,
\[
 |\mathcal A_m f(U\xi)-\mathcal A_m f(\xi)|
 \le \|f\circ U-f\|_\infty.
\]
Uniform continuity on the compact sphere and local rotations taking
one plane to a nearby plane prove continuity of $\mathcal A_m f$.
The identical argument on the Grassmannian proves continuity for the
adjoint. Positive lower bounds are preserved by averaging and by
positive powers. The family of planes through $u$ is the same as that
through $-u$, so the adjoint output is even.
\end{proof}

It is convenient to remove the fixed dimensional factor and write
\begin{equation}\label{eq:radial-P}
 \mathcal P_m(\rho)
 =\left\{\mathcal A_m^*[(\mathcal A_m(\rho^m))^m]\right\}^{1/(n-m)}.
\end{equation}
Then $\mathcal P_m(1)=1$ and
$\rho_{\mathbb I_mK}=C_{n,m}^{1/(n-m)}\mathcal P_m(\rho_K)$.
If $0<r\le\rho\le R$, positivity of averaging gives
\begin{equation}\label{eq:radial-bounds}
 r^{m^2/(n-m)}\le\mathcal P_m(\rho)\le R^{m^2/(n-m)}.
\end{equation}
The operator is orthogonally equivariant, maps star bodies to even star
bodies, and has homogeneity $m^2/(n-m)$. These statements use only
\eqref{eq:operator-normalized}, not a fixed-point hypothesis.

\subsection{The section moment as a radial functional}

For a positive continuous function $\rho$, let $K_\rho$ be its star body
and put
\[
 \begin{gathered}
s_\rho=\mathcal R_m(\rho^m/m),\qquad V(\rho)=\frac1n\int_{\Sph^{n-1}}\rho^n\,du,\\
 \Ph(\rho):=\Ph(K_\rho)=\frac1{m+1}\int_{G(n,m)}s_\rho^{m+1}\,d\nu_{n,m}.
\end{gathered}
\]
We use the same symbol $\Ph$ for the body functional and its radial
realization. If $K$ is merely a compact Borel set, the section formula
\eqref{intro:phi} still defines $\Ph(K)$. Section-volume measurability
follows by integrating the Borel indicator in local orthonormal frames
of the tautological bundle. If $K\subset RB_2^n$, then
$|K\cap\xi|_m\le\omega_mR^m$, so the functional is finite.

\begin{proposition}[Radial first variation and stationarity]
\label{prop:first-variation}
On the open cone of positive continuous functions, $\Ph$ and $V$ are
smooth in the supremum norm. Their derivatives are
\begin{align}
 D\Ph(\rho)[h]
 &=\int_{\Sph^{n-1}}\mathcal R_m^*(s_\rho^m)\rho^{m-1}h\,du,
 \label{eq:first-variation}\\
 DV(\rho)[h]&=\int_{\Sph^{n-1}}\rho^{n-1}h\,du.\notag
\end{align}
The positive stationary points of
$\Ph(\rho)^{1/(m+1)}/V(\rho)^{m/n}$ are precisely the radial
functions satisfying $\mathbb I_mK_\rho=cK_\rho$ for some $c>0$.
At such a point,
\begin{equation}\label{eq:scaling-stationarity}
 m(m+1)\Ph(K)=n c^{n-m}V(K).
\end{equation}
\end{proposition}
\begin{proof}
The map $\rho\mapsto\rho^m$ is a continuous polynomial on
$C(\Sph^{n-1})$. The transform $\mathcal R_m$ is bounded in the
supremum norm, with norm $\sigma_m$, and integration is a bounded
linear functional. The chain rule therefore gives
\[
 D\Ph(\rho)[h]
 =\int_{G(n,m)}s_\rho^m\mathcal R_m(\rho^{m-1}h)\,d\nu_{n,m}.
\]
Adjointness gives \eqref{eq:first-variation}. The formula for $DV$
follows in the same way.

The derivative of the positive quotient vanishes in every continuous
direction if and only if
\[
 D\Ph(\rho)=\lambda DV(\rho),\qquad
 \lambda=\frac{m(m+1)}{n}\frac{\Ph(\rho)}{V(\rho)}>0.
\]
By \eqref{eq:first-variation}, this means that the continuous function
$\mathcal R_m^*(s_\rho^m)\rho^{m-1}-\lambda\rho^{n-1}$ has zero
integral against every continuous $h$. Taking $h$ to be that function
shows it vanishes everywhere. Since $\rho>0$, division by
$\rho^{m-1}$ yields
$\mathcal R_m^*(s_\rho^m)=\lambda\rho^{n-m}$.
This is the fixed-point equation with $c=\lambda^{1/(n-m)}$.
Conversely, that equation gives $D\Ph=c^{n-m}DV$.
Evaluation in the scaling direction $h=\rho$ proves
\eqref{eq:scaling-stationarity} and identifies the required multiplier.
\end{proof}

For convex $K$ the one-way measure identity \eqref{intro:cone-relation}
follows by substituting the fixed-point equation into the definition of
$\widetilde A_{m,m+1}$ in \cite[Definition 3.1]{CLWX2026}.
No inversion of the reverse radial Gauss image is used in this argument.
The classification below concerns the pointwise radial equation, not an
unproved converse to that measure identity.

\subsection{A uniform quadratic remainder}

The next estimate is formulated for arbitrary small radial increments.
It will apply to a volume-preserving path even when its time derivative
has not been established.

\begin{proposition}[Quadratic remainder at a fixed point]
\label{prop:quadratic-stationarity}
Let $K$ be any star body satisfying $\mathbb I_mK=cK$ with $c>0$,
and write $\rho=\rho_K$. There exist $\varepsilon_1,C>0$ such that
for every continuous $h$ with $\|h\|_\infty\le\varepsilon_1$,
\begin{equation}\label{eq:quadratic-stationarity}
 \left|\Ph(\rho+h)-\Ph(\rho)
 -c^{n-m}\bigl(V(\rho+h)-V(\rho)\bigr)\right|
 \le C\|h\|_\infty^2.
\end{equation}
In particular, if a family $K_t$ has $|K_t|=|K|$ and
$\|\rho_{K_t}-\rho_K\|_\infty\le Mt$, then
\begin{equation}\label{eq:quadratic-path}
 |\Ph(K_t)-\Ph(K)|\le CM^2t^2
\end{equation}
for small $t$, without any assumption that $\partial_t\rho_{K_t}$ exists.
\end{proposition}
\begin{proof}
We give an explicit uniform remainder, rather than write a formal first
variation along a possibly nondifferentiable path. Let
\[
 \rho_{\min}=\min\rho>0,\quad R_0=\max\rho,\quad
 \varepsilon_1=\min\{1,\rho_{\min}/2\},\quad B=R_0+\varepsilon_1,
 \quad S=\omega_mB^m.
\]
For a fixed continuous $h$ with $H:=\|h\|_\infty\le\varepsilon_1$, introduce
an auxiliary scalar $\tau\in[0,1]$ and set
\[
 r_\tau=\rho+\tau h,\qquad s_\tau=\mathcal R_m(r_\tau^m/m).
\]
Then $\rho_{\min}/2\le r_\tau\le B$ and $0<s_\tau\le S$. These bounds are
uniform in $h$ and $\tau$. All powers are integer powers, so their
derivatives in $\tau$ exist as derivatives of continuous functions in
supremum norm. In particular,
\[
 s_\tau'=\mathcal R_m(r_\tau^{m-1}h),\qquad
 \|s_\tau'\|_\infty\le\sigma_mB^{m-1}H.
\]
For $m\ge2$,
\[
 s_\tau''=(m-1)\mathcal R_m(r_\tau^{m-2}h^2),\qquad
 \|s_\tau''\|_\infty\le(m-1)\sigma_mB^{m-2}H^2;
\]
when $m=1$ this second derivative is identically zero.

Define the ordinary scalar function
\[
 Q(\tau)=\Ph(r_\tau)-c^{n-m}V(r_\tau).
\]
Uniform bounds justify differentiation under the finite-measure integrals.
At $\tau=0$, adjointness and the fixed-point equation give
\begin{align*}
 Q'(0)
 &=\int_{G(n,m)}s_\rho^m\mathcal R_m(\rho^{m-1}h)\,d\nu
       -c^{n-m}\int_{\Sph^{n-1}}\rho^{n-1}h\,d\HH^{n-1}\\
 &=\int_{\Sph^{n-1}}
   [\mathcal R_m^*(s_\rho^m)\rho^{m-1}-c^{n-m}\rho^{n-1}]
        h\,d\HH^{n-1}=0.
\end{align*}
The last equality is pointwise: the exponent is
$(n-m)+(m-1)=n-1$. This is the stationarity identity of Proposition~\ref{prop:first-variation}.

For the second derivative,
\begin{align*}
 \frac{d^2}{d\tau^2}\Ph(r_\tau)
 &=m\int_{G(n,m)}s_\tau^{m-1}(s_\tau')^2\,d\nu
       +\int_{G(n,m)}s_\tau^m s_\tau''\,d\nu,\\
 \frac{d^2}{d\tau^2}V(r_\tau)
 &=(n-1)\int_{\Sph^{n-1}}r_\tau^{n-2}h^2\,d\HH^{n-1}.
\end{align*}
Let
\[
 C_{\Phi}=mS^{m-1}\sigma_m^2B^{2m-2}
 +\begin{cases}(m-1)S^m\sigma_mB^{m-2},&m\ge2,\\0,&m=1,\end{cases}
 \qquad C_V=(n-1)\sigma_nB^{n-2}.
\]
Since $\nu$ is a probability measure,
$|Q''(\tau)|\le(C_{\Phi}+c^{n-m}C_V)H^2$ on $[0,1]$.
The one-variable integral form of Taylor's theorem now gives
\[
 Q(1)-Q(0)=Q'(0)+\int_0^1(1-\tau)Q''(\tau)\,d\tau,
\]
so the claimed inequality holds with
$C=(C_{\Phi}+c^{n-m}C_V)/2$.

Finally put $h=\rho_{K_t}-\rho$. For small $t$ its supremum norm is
at most $\varepsilon_1$, while volume preservation makes
$V(\rho+h)-V(\rho)=0$. The result follows from $H\le Mt$.
The parameter $\tau$ used in this proof is a straight-line interpolation
of functions for a fixed $h$. It is not the Steiner time $t$, and nothing
here assumes existence of a derivative of $t\mapsto\rho_{K_t}$.
\end{proof}

\subsection{The line-section case and the critical endpoint}

\begin{proposition}[The $m=1$ star-body fixed-point equation]\label{prop:m1}
If $n\ge3$ and $\mathbb I_1K=cK$ for a star body $K$ and $c>0$,
then $K$ is an origin-centered Euclidean ball.
\end{proposition}
\begin{proof}
A line section has length $\rho(u)+\rho(-u)$, and
$\mathcal R_1^*=(2/\sigma_n)\mathcal A_1^*$. Thus
\[
 c^{n-1}\rho(u)^{n-1}
 =\frac2{\sigma_n}(\rho(u)+\rho(-u)).
\]
The right side is even, so strict positivity first gives
$\rho(u)=\rho(-u)$ and then
$c^{n-1}\rho(u)^{n-2}=4/\sigma_n$. Hence $\rho$ is constant.
\end{proof}

For comparison with the classical operator, let
$\rho_{IK}(v)=|K\cap v^\perp|_{n-1}$. A hyperplane containing $u$
is $v^\perp$ with $v\perp u$. Since the integrand is even,
\[
 \rho_{\mathbb I_{n-1}K}(u)
 =\frac1{\sigma_n}\int_{\Sph^{n-1}\cap u^\perp}
     \rho_{IK}(v)^{n-1}\,d\mathcal H^{n-2}(v)
 =\frac{n-1}{\sigma_n}\rho_{I^2K}(u).
\]
This proves \eqref{intro:endpoint}, including its normalization. The
ellipsoidal endpoint classification in \cite[Theorem 1.1]{MSY2025}
explains why the ball theorem is restricted to $m\le n-2$.

\section{Radon regularity and smoothness of fixed points}\label{sec:regularity}

In this section $n\ge3$ and $2\le m<n$. We work with the probability
averages of Section~\ref{sec:preliminaries}. The main objective is to
remove any initial Lipschitz assumption from the fixed-point problem.
The proof uses the normal operator, then full-domain Sobolev estimates,
and finally a nonlinear bootstrap with separate supremum bounds.

\subsection{The positive normal kernel}

\begin{lemma}[An explicit positive normal kernel]\label{lem:normal-kernel}
Let $n\ge3$ and $2\le m\le n-1$. Put $N_m=\mathcal A_m^*\mathcal A_m$
and
\[
 \vartheta_{n,m}
 =\frac{\Gamma(m/2)\Gamma((n-1)/2)}
 {\Gamma((m-1)/2)\Gamma(n/2)}>0.
\]
For bounded measurable $f$, almost everywhere in $u$,
\begin{equation}\label{eq:normal-kernel}
 N_m f(u)=\vartheta_{n,m}\int_{\Sph^{n-1}}
 \frac{f(v)}{(1-\langle u,v\rangle^2)^{(n-m)/2}}\,\dd\sigma(v).
\end{equation}
For continuous $f$ the identity holds everywhere. Its kernel is
integrable and has total mass one. In particular, for $f\ge0$,
\begin{equation}\label{eq:normal-positive}
 N_m f(u)\ge \vartheta_{n,m}\int f\,\dd\sigma.
\end{equation}
For the same nonnegative $f$,
\begin{equation}\label{eq:nonlinear-positive}
 \mathcal A_m^*[(\mathcal A_m f)^m](u)
 \ge (N_m f(u))^m
 \ge \left(\vartheta_{n,m}\int f\,\dd\sigma\right)^m.
\end{equation}
\end{lemma}
\begin{proof}
Fix $u$. First choose an $m$-plane $\xi$ through $u$ uniformly, and then
choose $v$ uniformly in $\Sph^{n-1}\cap\xi$. Write
\[
 v=t u+\sqrt{1-t^2}\,w,\qquad w\in\Sph^{n-1}\cap u^\perp.
\]
The coordinate $t$ has density
\[
 c_m(1-t^2)^{(m-3)/2},\quad -1<t<1,\qquad
 c_j=\frac{\Gamma(j/2)}{\sqrt\pi\Gamma((j-1)/2)}.
\]
After averaging over $\xi$, the conditional distribution of $w$ is
uniform by rotational invariance. Sphere probability itself has the
same conditional distribution of $w$ and coordinate density
$c_n(1-t^2)^{(n-3)/2}$. The ratio of these densities is precisely
$\vartheta_{n,m}(1-t^2)^{-(n-m)/2}$. This proves
\eqref{eq:normal-kernel}. Its endpoint integrability is equivalent to
$(m-3)/2>-1$, which holds because $m\ge2$.
Since $0<1-t^2\le1$ away from two measure-zero poles, the kernel is
at least $\vartheta_{n,m}$, giving \eqref{eq:normal-positive}.
Jensen's inequality for the conditional probability defining
$\mathcal A_m^*$, applied to the convex map $z\mapsto z^m$ on
$[0,\infty)$, proves \eqref{eq:nonlinear-positive}.
\end{proof}

\subsection{Spectrum and Sobolev estimates}

\begin{lemma}[Exact eigenvalues]\label{lem:eig}
For $n\ge3$, $2\le m\le n-1$, the operator $N_m$ annihilates odd
spherical harmonics. On $\HH_{2k}(\Sph^{n-1})$ it acts by
\begin{equation}\label{eq:eig}
 N_mY_{2k}=\lambda_{2k}Y_{2k},\qquad
 \lambda_{2k}
 =\frac{(\frac12)_k(\frac{n-m}{2})_k}
 {(\frac m2)_k(\frac{n-1}{2})_k},\quad k\ge0.
\end{equation}
Thus $\lambda_0=1$ and
\begin{equation}\label{eq:eig-asymptotic}
 c(1+k)^{-(m-1)}\le\lambda_{2k}
 \le C(1+k)^{-(m-1)}.
\end{equation}
\end{lemma}
\begin{proof}
\textbf{1. Why it suffices to compute one scalar per degree.}
The space $\HH_k$ consists of the restrictions to the sphere of homogeneous
harmonic polynomials of degree $k$. For $n\ge3$ these are the irreducible,
pairwise inequivalent spherical harmonic representations of the orthogonal
group. The averaging operators commute with rotations, so $N_m$ acts as a
scalar on each $\HH_k$. Equivalently, this follows from the Funk--Hecke
principle for the zonal kernel of Lemma~\ref{lem:normal-kernel}.
If $k$ is odd, every section sphere contains $v$ and $-v$ with equal measure,
and the section average vanishes. Thus only degrees $2k$ remain.

\textbf{2. Use a normalized zonal harmonic.}
Fix $u\in\Sph^{n-1}$ and set $\zeta=(n-2)/2>0$. The function
\[
 Y(v)=P_{2k}(u\cdot v),\qquad
 P_{2k}(t)=\frac{C_{2k}^{\zeta}(t)}{C_{2k}^{\zeta}(1)}
\]
is a spherical harmonic of degree $2k$ and satisfies $Y(u)=1$.
Consequently its eigenvalue is $(N_mY)(u)$. By the normal-kernel
disintegration this is $\mathbb E P_{2k}(T)$, where
$T$ has density $c_m(1-t^2)^{(m-3)/2}$ on $(-1,1)$.
Writing $(b)_j=b(b+1)\cdots(b+j-1)$ and $(b)_0=1$, the beta integral gives
\[
 \mathbb ET^{2j}
 =\frac{B(j+1/2,(m-1)/2)}{B(1/2,(m-1)/2)}
 =\frac{(1/2)_j}{(m/2)_j}.
\]

\textbf{3. Evaluate a finite sum, not an infinite interchange.}
The even Gegenbauer polynomial has the terminating expansion
\[
 C_{2k}^{\zeta}(t)
 =(-1)^k\frac{(\zeta)_k}{k!}
   \sum_{j=0}^{k}\frac{(-k)_j(k+\zeta)_j}{(1/2)_j\,j!}\,t^{2j}.
\]
This formula is finite, so term-by-term averaging has no convergence
requirement. The preceding moment identity cancels $(1/2)_j$ and yields
\[
 \mathbb E C_{2k}^{\zeta}(T)
 =(-1)^k\frac{(\zeta)_k}{k!}
   \sum_{j=0}^{k}\frac{(-k)_j(k+\zeta)_j}{(m/2)_j\,j!}.
\]
The finite Chu--Vandermonde identity
$\sum_{j=0}^k(-k)_j(b)_j/((c)_j j!)=(c-b)_k/(c)_k$
with $b=k+\zeta$, $c=m/2>0$ reduces this to
\[
 (-1)^k\frac{(\zeta)_k}{k!}
    \frac{(m/2-k-\zeta)_k}{(m/2)_k}
 =\frac{(\zeta)_k}{k!}\frac{((n-m)/2)_k}{(m/2)_k}.
\]
For the last equality reverse the order of the $k$ factors in the
numerator; each contributes a minus sign.
To recall why the finite summation identity holds for these parameters,
first take $0<b<c$ and average the binomial expansion of $(1-T)^k$
against a $\mathrm{Beta}(b,c-b)$ random variable. This gives
$(c-b)_k/(c)_k$. For fixed $c>0$, both sides are polynomials in $b$,
so the identity extends to every real $b$. Thus the substitution
$b=k+\zeta$ is valid even when $b\ge c$.

\textbf{4. Divide by the normalization and estimate.}
The duplication formula for Pochhammer products gives
\[
 C_{2k}^{\zeta}(1)=\frac{(2\zeta)_{2k}}{(2k)!}
 =\frac{(\zeta)_k(\zeta+1/2)_k}{(1/2)_k k!}.
\]
Dividing the expectation by this expression proves \eqref{eq:eig},
since $\zeta+1/2=(n-1)/2$. All factors are positive, and $k=0$ gives $1$.
Finally $(b)_k=\Gamma(k+b)/\Gamma(b)$ and the Gamma-ratio asymptotic
give the two exponents
\[
 \tfrac12-\tfrac m2=\tfrac{1-m}{2},\qquad
 \tfrac{n-m}{2}-\tfrac{n-1}{2}=\tfrac{1-m}{2}.
\]
Their sum is $1-m$. This proves the two-sided estimate after adjusting
the constants to include finitely many small values of $k$.
\end{proof}

To avoid a metric normalization ambiguity, define the elliptic operator
on either compact manifold $M=\Sph^{n-1}$ or $\GG(n,m)$ by
\[
 \mathcal L_M=-\sum_{1\le i<j\le n}(X_{ij}^M)^2,
\]
where $X_{ij}^M$ is the infinitesimal action of unit-speed rotation in
the $(i,j)$ coordinate plane. These fields span the tangent space at
every point and are skew-adjoint on smooth functions for invariant
measure. We take the nonnegative self-adjoint Friedrichs realization
of this elliptic differential operator. Then
\[
 \|h\|_{H^s(M)}=\|(1+\mathcal L_M)^{s/2}h\|_2
\]
is a Sobolev norm equivalent to the local-chart norm; see, for the
homogeneous-space setting, \cite{Helgason1984}. On the sphere $\mathcal L_M=-\Delta_{\Sph}$.
This convention fixes the exact norm used later for the local contraction.
The averaging operators send smooth functions to smooth functions:
in a local smooth choice of rotations, derivatives fall on the smooth
input inside an integral over a fixed compact fiber. Equivariance then
gives $X_{ij}^{\GG}\mathcal A_m=\mathcal A_mX_{ij}^{\Sph}$ and hence
$\mathcal L_{\GG}\mathcal A_m=\mathcal A_m\mathcal L_{\Sph}$ on smooth inputs.

Set $a=(m-1)/2$. The frequency weight for a gain of $a$ Sobolev
orders is $(1+\kappa_{2k})^a$, where $\kappa_j=j(j+n-2)$.
By \eqref{eq:eig-asymptotic}, the product
$(1+\kappa_{2k})^a\lambda_{2k}$ is uniformly bounded. This is the
precise balance used in both estimates below.

\begin{proposition}[Sobolev gains on both full domains]
\label{prop:sobolev-radon}
Let $n\ge3$, $2\le m\le n-1$ and $a=(m-1)/2$. For every $s\ge0$,
\begin{equation}\label{eq:sobolev-gains}
 \begin{split}
 \|\mathcal A_mf\|_{H^{s+a}(\GG(n,m))}
 &\le C\|f\|_{H^s(\Sph^{n-1})},\\
 \|\mathcal A_m^*g\|_{H^{s+a}(\Sph^{n-1})}
 &\le C\|g\|_{H^s(\GG(n,m))}.
 \end{split}
\end{equation}
The second estimate holds for arbitrary $g$ in its displayed domain,
not only for functions in the range of $\mathcal A_m$.
\end{proposition}
\begin{proof}
Put $\kappa_k=k(k+n-2)$. Choose an orthonormal basis
$Y_{2k,j}$ of $\HH_{2k}$ in sphere probability and set
\[
 Z_{k,j}=\lambda_{2k}^{-1/2}\mathcal A_mY_{2k,j}.
\]
Adjointness and \eqref{eq:eig} show that the $Z_{k,j}$ form an
orthonormal family on $\GG(n,m)$. Each satisfies
$\mathcal L_{\GG}Z_{k,j}=\kappa_{2k}Z_{k,j}$.
These vectors span the closure of the range of $\mathcal A_m$;
there is no assertion that they exhaust all Grassmannian eigenfunctions.

For $f=\sum f_k$ the first norm is exactly
\[
 \|\mathcal A_m f\|_{H^{s+a}}^2
 =\sum_{k\ge0}(1+\kappa_{2k})^{s+a}
          \lambda_{2k}\|f_{2k}\|_2^2.
\]
Since $(1+\kappa_{2k})^a\lambda_{2k}$ is uniformly bounded by
\eqref{eq:eig-asymptotic}, this proves the first estimate.
For arbitrary $g$, adjointness gives the exact coefficient identity
\[
 \langle\mathcal A_m^*g,Y_{2k,j}\rangle
 =\lambda_{2k}^{1/2}\langle g,Z_{k,j}\rangle.
\]
Every odd coefficient vanishes, and any part of $g$ orthogonal to
all the $Z_{k,j}$ is annihilated. Therefore
\begin{align*}
 \|\mathcal A_m^*g\|_{H^{s+a}}^2
 &=\sum_{k,j}(1+\kappa_{2k})^{s+a}
       \lambda_{2k}|\langle g,Z_{k,j}\rangle|^2\\
 &\le C\sum_{k,j}(1+\kappa_{2k})^s
              |\langle g,Z_{k,j}\rangle|^2
 \le C\|g\|_{H^s(\GG)}^2.
\end{align*}
For completeness, put $G=(1+\mathcal L_{\GG})^{s/2}g$.
Because $Z_{k,j}$ is an eigenfunction of $\mathcal L_{\GG}$,
\[
 (1+\kappa_{2k})^{s/2}\langle g,Z_{k,j}\rangle
       =\langle G,Z_{k,j}\rangle.
\]
The final sum is therefore at most $\|G\|_2^2$ by ordinary Bessel
inequality. The family $Z_{k,j}$ need not be a basis of all Grassmannian
functions: any orthogonal remainder is killed by the adjoint. This is
why the estimate applies after a nonlinear power, even if that power
has left the range of $\mathcal A_m$.
Smooth approximation extends the calculation to all the stated inputs.
\end{proof}

\subsection{Bounded products and the nonlinear bootstrap}

\begin{lemma}[A Sobolev algebra with an independent supremum norm]
\label{lem:powers}
Let $M$ be a compact smooth manifold without boundary and $s\ge0$.
Then
\[
 \mathfrak B^s(M)=H^s(M)\cap L^\infty(M),\qquad
 \|f\|_{\mathfrak B^s}=\|f\|_{H^s}+\|f\|_\infty
\]
is a Banach algebra. In particular, for every positive integer $j$,
$f\mapsto f^j$ is a continuous polynomial, hence a smooth map, on
$\mathfrak B^s(M)$. We have the tame product estimate
\begin{equation}\label{eq:tame-product}
 \|fg\|_{H^s}\le C_s
 \bigl(\|f\|_\infty\|g\|_{H^s}
       +\|g\|_\infty\|f\|_{H^s}\bigr).
\end{equation}
If $k\ge0$ is an integer, $f\in H^k(M)\cap L^\infty(M)$ takes
values in a compact interval $I$, and $F$ is smooth on a neighborhood
of $I$, then $F(f)\in H^k(M)$ and
\begin{equation}\label{eq:integer-composition}
 \|F(f)\|_{H^k}\le C_{M,k,F,I}(1+\|f\|_{H^k}).
\end{equation}
\end{lemma}
\begin{proof}
For $s=0$, \eqref{eq:tame-product} is immediate. For $s>0$ start
with the Euclidean inhomogeneous Kato--Ponce estimate
\cite[Theorem 1, inequality (3)]{GrafakosOh2014}, with target exponent $2$
and the two exponent pairs $(\infty,2)$ and $(2,\infty)$.
To pass to $M$, choose a finite smooth atlas, a subordinate
partition of unity $\chi_\alpha$, and cutoffs $\psi_\alpha$ supported
in the same charts and equal to one near $\operatorname{supp}\chi_\alpha$.
In a chart $\kappa_\alpha$, extend
$f_\alpha=(\chi_\alpha f)\circ\kappa_\alpha^{-1}$ and
$g_\alpha=(\psi_\alpha g)\circ\kappa_\alpha^{-1}$ by zero.
They have compact support inside the chart, and
\[
 (\chi_\alpha fg)\circ\kappa_\alpha^{-1}=f_\alpha g_\alpha.
\]
Applying the Euclidean estimate, summing over finitely many charts,
and using the equivalence of the local-chart Sobolev norm and the
spectral norm of a positive elliptic operator proves
\eqref{eq:tame-product}. Multiplication by the fixed cutoffs is bounded
on $H^s$; their supremum norms give the required independent
$L^\infty$ factors. Adding the elementary supremum product bound
makes multiplication continuous on $\mathfrak B^s$. Completeness
follows because an $H^s$-Cauchy and $L^\infty$-Cauchy sequence has
the same $L^2$ limit on this finite-volume manifold.
For nonsmooth bounded inputs, heat regularization preserves a uniform
supremum bound, converges in $H^s$, and allows passage to the product
by $L^2$ convergence and weak $H^s$ compactness. Density in the
$L^\infty$ norm is neither asserted nor needed.

Only integer-order smooth composition is needed in the bootstrap.
For smooth $f$ and integer $k\ge1$, the Gagliardo--Nirenberg
inequalities on Euclidean charts, with cutoffs as above and the
lower-order compact-manifold term included, give for $1\le j\le k$,
\[
 \|\nabla^j f\|_{L^{2k/j}}
 \le C\bigl(\|f\|_\infty^{1-j/k}\|f\|_{H^k}^{j/k}
                    +\|f\|_\infty\bigr).
\]
A derivative of order $\ell\le k$ of $F(f)$ is a sum of terms
$F^{(r)}(f)\prod_{i=1}^r\nabla^{j_i}f$, where $j_i\ge1$ and
$\sum j_i=\ell$. The derivatives of $F$ are bounded on $I$.
H\"older with exponents $2k/j_i$, followed by the finite-volume
inclusion $L^{2k/\ell}\subset L^2$, bounds each term by
$C(1+\|f\|_{H^k})$. This proves \eqref{eq:integer-composition}
for smooth inputs. Heat regularization of $f$ stays in $I$ by the
maximum principle and converges in $H^k$. The estimates just proved
are uniform, and weak compactness together with $L^2$ convergence
of $F(f)$ proves the assertion in general. For $k=0$, boundedness
of $F$ on $I$ suffices.
\end{proof}

\begin{remark}[Why no Grassmannian Sobolev embedding is being used]
Although $H^s(\Sph^{n-1})\subset L^\infty$ when $s>(n-1)/2$,
one cannot infer $H^{s+a}(\GG(n,m))\subset L^\infty$ from this:
$\dim\GG(n,m)=m(n-m)$ can be much larger than $n-1$.
For example, $n=6,m=3,s=3$ gives $s+a=4<9/2$.
Instead the independent bound
$\|\mathcal A_mf\|_\infty\le\|f\|_\infty$ supplies the
supremum norm in $\mathfrak B^{s+a}(\GG)$.
\end{remark}

Recall $\mathcal P_m$ from \eqref{eq:radial-P}. Its two averaging
steps improve regularity, while its powers preserve the independent
supremum bounds. The final root is applied only after a positive lower
bound has been established.

\begin{proposition}[One nonlinear step adds $m-1$ Sobolev orders]
\label{prop:nonlinear-gain}
Let $2\le m<n$, let $k\ge0$ be an integer, and suppose
$\rho\in H^k(\Sph^{n-1})$ with $0<r\le\rho\le R<\infty$ almost
everywhere. Then
\[
 \mathcal P_m(\rho)\in H^{k+m-1}(\Sph^{n-1})\cap L^\infty,
 \qquad r^{m^2/(n-m)}\le\mathcal P_m(\rho)\le R^{m^2/(n-m)}.
\]
\end{proposition}
\begin{proof}
Put $a=(m-1)/2$. The precise sequence is
\[
 \begin{aligned}
 \rho^m&\in\mathfrak B^k(\Sph^{n-1}),\\
 v=\mathcal A_m(\rho^m)&\in\mathfrak B^{k+a}(\GG(n,m)),\\
 v^m&\in\mathfrak B^{k+a}(\GG(n,m)),\\
 w=\mathcal A_m^*(v^m)&\in H^{k+2a}(\Sph^{n-1})\cap L^\infty.
 \end{aligned}
\]
At every stage the supremum bound comes from averaging or an integer
power; it does not come from a dimension-dependent embedding.
Also $r^{m^2}\le w\le R^{m^2}$. Since $k+2a=k+m-1$ is an integer,
apply \eqref{eq:integer-composition} to
$F(t)=t^{1/(n-m)}$ on $[r^{m^2},R^{m^2}]$.
This proves the regularity gain and the pointwise bounds.
\end{proof}

\begin{theorem}[Automatic positivity and smoothness of bounded solutions]
\label{thm:smooth}
Let $n\ge3$, $2\le m\le n-1$, and let
$0\le\rho\in L^\infty(\Sph^{n-1})$ be nonzero as an
almost-everywhere class. Suppose, for some $\Lambda>0$,
\begin{equation}\label{eq:bounded-equation}
 \mathcal A_m^*[(\mathcal A_m(\rho^m))^m]
 =\Lambda\rho^{n-m}\qquad\text{almost everywhere}.
\end{equation}
Then $\rho$ has an even, strictly positive $C^\infty$ representative.
In particular, every star body satisfying $\mathbb I_mK=cK$ is smooth.
\end{theorem}
\begin{proof}
Apply \eqref{eq:nonlinear-positive} to $f=\rho^m$.
Since $\int\rho^m\,\dd\sigma>0$, equation
\eqref{eq:bounded-equation} yields the uniform lower bound
\begin{equation}\label{eq:rho-positive}
 \rho(u)\ge r_*:=\left[
 \frac{(\vartheta_{n,m}\int\rho^m\,\dd\sigma)^m}{\Lambda}
 \right]^{1/(n-m)}>0
 \quad\text{almost everywhere}.
\end{equation}
The assumed essential upper bound and this positive lower bound
supply a fixed compact range for every root composition.
Now $\rho\in H^0$. If $\rho\in H^k$ for an integer $k\ge0$,
Proposition~\ref{prop:nonlinear-gain} and
$\mathcal P_m(\rho)=\Lambda^{1/(n-m)}\rho$ imply
$\rho\in H^{k+m-1}$. Inductively,
\begin{equation}\label{eq:rho-smooth}
 \rho\in H^{j(m-1)}(\Sph^{n-1})\quad(j=0,1,2,\ldots).
\end{equation}
For every integer $\ell\ge0$, choose $j$ with
$j(m-1)>\ell+(n-1)/2$. Sobolev embedding gives a $C^\ell$
representative. Choose one continuous representative once $j(m-1)>(n-1)/2$.
For every $\ell$ a higher-order representative agrees with it almost
everywhere. Two continuous functions agreeing almost everywhere agree
everywhere, so the same chosen representative is $C^\ell$ for every
$\ell$. It is therefore $C^\infty$. It is bounded below by $r_*$
everywhere: a strict violation at one point would, by continuity,
give a positive-measure violation of the almost-everywhere lower bound. The
right-hand side defining $\mathcal A_m^*$ is even, hence so is
$\rho$; the now-continuous equation holds everywhere by
Lemma~\ref{lem:averaging}.

For a star-body fixed point, \eqref{eq:operator-normalized} reads
\begin{equation}\label{eq:regularity-fp}
 c^{n-m}\rho^{n-m}
 =C_{n,m}\mathcal A_m^*[(\mathcal A_m(\rho^m))^m].
\end{equation}
Take $\Lambda=c^{n-m}/C_{n,m}$ in the preceding argument. Because
its original radial function is continuous, it already coincides
everywhere with the smooth representative.
\end{proof}

\begin{corollary}[Finite-iterate regularization]\label{cor:iterate-smooth}
Let $2\le m<n$ and let $\rho_0>0$ be continuous. Set
$\rho_{j+1}=\mathcal P_m(\rho_j)$. Then
$\rho_j\in H^{j(m-1)}$ for every $j\ge0$. In particular,
$\rho_j\in C^\ell$ whenever $j(m-1)>\ell+(n-1)/2$.
\end{corollary}
\begin{proof}
Each fixed iterate is continuous, bounded and bounded away from zero
by Lemma~\ref{lem:averaging}. Apply
Proposition~\ref{prop:nonlinear-gain} inductively and then Sobolev embedding.
No bounds uniform in $j$ are asserted here.
\end{proof}

\begin{remark}
The bounded-solution conclusion requires $m\ge2$. For $m=1$, an
antipodally symmetric measurable set $D$ and a function $b\one_D$
can satisfy the corresponding equation after choosing its scalar
coefficient, without being strictly positive or continuous.
The star-body case $m=1$ is covered separately by
Proposition~\ref{prop:m1}. Moreover, finite-iterate regularization gives
no estimates uniform in the number of iterates, and therefore does not
imply compactness or convergence of arbitrary normalized orbits.
\end{remark}

\section{Graph coordinates and weighted section representations}\label{sec:graph}\subsection{Haar measure in graph coordinates}

Fix $u\in\Sph^{n-1}$. For $E\in\GG(u^\perp,m)$ and $a\in E$, set
\begin{equation}\label{eq:graph}
 \xi(E,a)=\bigl\{y+\langle a,y\rangle u:y\in E\bigr\}
 \in\GG(n,m).
\end{equation}
Every $m$-plane not containing $u$ has a unique representation of this
form. The excluded set of planes has Haar measure zero.

The plane $E$ is the projection of $\xi$ onto $u^\perp$, and $a$ records its slope: at horizontal point $y$ the
height of $\xi$ is $\langle a,y\rangle$. The projection is one-to-one exactly
when $\xi$ does not contain $u$. Thus its inverse is linear, and its vertical
component is a linear functional on $E$, uniquely represented by $a$.
The formula below computes the distribution of both $E$ and this slope;
it is not obtained by assuming that two marginal distributions are independent.

\begin{lemma}[Haar probability in graph coordinates]\label{lem:Haar}
For every nonnegative measurable $F$ on $\GG(n,m)$,
\begin{equation}\label{eq:haar}
 \int_{\GG(n,m)}F(\xi)\,\dd\nu_{n,m}(\xi)
 =A_{n,m}\int_{\GG(u^\perp,m)}\int_E
 F(\xi(E,a))\,(1+|a|^2)^{-n/2}\,\dd a\,
 \dd\nu_{u^\perp,m}(E),
\end{equation}
where
\begin{equation}\label{eq:A}
 A_{n,m}=\frac{\Gamma(n/2)}{\pi^{m/2}\Gamma((n-m)/2)}.
\end{equation}
\end{lemma}
\begin{proof}
Let $d=n-1$, let $H=O(u^\perp)$ with Haar probability, and let
$\eta$ be the pushforward of $\nu_{n,m}$ to the graph bundle
$\{(E,a):E\in\GG(u^\perp,m),\ a\in E\}$. The omitted set has
measure zero: for a Haar plane, $Z=|P_\xi u|^2$ has law
$\mathrm{Beta}(m/2,(n-m)/2)$ and has no atom at $1$.
To see this law, rotate a fixed $m$-plane by a Haar orthogonal matrix.
Then $Z$ is the sum of the first $m$ squared coordinates of a uniform
point on $\Sph^{n-1}$, equivalently the ratio of independent chi-square
variables of degrees $m$ and $n-m$ to their sum.

The graph coordinates are equivariant for $H$. For every nonnegative
Borel function $G$ on this bundle, invariance and Tonelli give
\begin{equation}\label{eq:haar-orbit-average}
 \int G\,\dd\eta
 =\int\!\left[\int_H G(UE,Ua)\,\dd U\right]\dd\eta(E,a).
\end{equation}
For each fixed $r>0$, $H$ acts transitively on
\[
 \mathcal O_r=\{(E,a):a\in E,\ |a|=r\}.
\]
Indeed, send $a/r$ to any prescribed unit vector and extend orthonormal
bases of the two $m$-planes and their orthogonal complements.
The orbit average in brackets therefore depends only on $r=|a|$.
Its value is
\[
 M_G(r)=\int_{\GG(u^\perp,m)}
       \int_{\Sph(E)}G(E,r\theta)\,\dd\sigma_E(\theta)\,\dd\nu(E),
\]
where $\sigma_E$ is probability on the unit sphere of $E$.
To justify this identification, the displayed probability is
$H$-invariant on $\mathcal O_r$; on a compact transitive $H$-space,
any invariant probability integrates a continuous function to its
Haar orbit average, hence is unique. Extension to nonnegative Borel
functions follows by the uniqueness of the resulting measures.
Thus the joint distribution is determined by its orbit average and
radial law, with no independence assumption.

The graph map has Gram matrix $I_E+a\otimes a$. Hence
\[
 Z=a^T(I_E+a\otimes a)^{-1}a=\frac{|a|^2}{1+|a|^2}.
\]
More explicitly, the beta density is
\[
 \frac{z^{m/2-1}(1-z)^{(n-m)/2-1}}{B(m/2,(n-m)/2)}\,dz,
 \qquad 0<z<1.
\]
Substitute $z=r^2/(1+r^2)$, so $dz=2r(1+r^2)^{-2}dr$.
The numerator contributes $r^{m-2}$, and all factors of $1+r^2$
together have exponent $-n/2$. This gives
\[
 \dd\mathbb P_{|a|}(r)
 =\frac{2}{B(m/2,(n-m)/2)}
   r^{m-1}(1+r^2)^{-n/2}\,\dd r.
\]
There is no atom at $r=0$. Insert this density and $M_G(r)$ into
\eqref{eq:haar-orbit-average}, and use polar coordinates
$\dd a=\sigma_mr^{m-1}\dd r\dd\sigma_E$. The coefficient is
\[
 \frac{2}{\sigma_m B(m/2,(n-m)/2)}
 =\frac{\Gamma(n/2)}{\pi^{m/2}\Gamma((n-m)/2)}.
\]
Finally choose $G(E,a)=F(\xi(E,a))$ to obtain \eqref{eq:haar}.
\end{proof}

\begin{theorem}[Weighted graph representation]\label{thm:graph}
For every compact Borel set $K\subset\R^n$ and
$1\le m\le n-1$, set $\beta=(n-m-1)/2\ge0$. Then
\begin{equation}\label{eq:weighted-graph}
 \begin{aligned}
 \Ph(K)=\frac{A_{n,m}}{m+1}
 \int_{\GG(u^\perp,m)}\int_E &(1+|a|^2)^{-\beta}\\
 &\times\left[\int_E\one_K(y+\langle a,y\rangle u)\,\dd y\right]^{m+1}
 \dd a\,\dd\nu_{u^\perp,m}(E).
 \end{aligned}
\end{equation}
\end{theorem}
\begin{proof}
Fix $(E,a)$ and let $T_a:E\to\R^n$ be $T_a y=y+\langle a,y\rangle u$.
For an orthonormal basis $e_1,\ldots,e_m$ of $E$,
\[
 \langle T_a e_i,T_a e_j\rangle=\delta_{ij}+a_i a_j.
\]
Thus the Gram matrix is $I+a\otimes a$. Its eigenvalues are $1$ in
directions perpendicular to $a$ and $1+|a|^2$ along $a$; the same determinant
formula holds at $a=0$. The area multiplier is therefore
\[
 J_mT_a=\sqrt{\det(I+a\otimes a)}=\sqrt{1+|a|^2}.
\]
Writing
\[
 Q_K(E,a)=\int_E\one_K(y+\langle a,y\rangle u)\,dy,
\]
the ordinary area formula for this injective linear map gives
$|K\cap\xi(E,a)|_m=\sqrt{1+|a|^2}\,Q_K(E,a)$.

Now insert this into the defining integral for $\Ph$ and apply
Lemma~\ref{lem:Haar}. The section-volume power contributes
$(1+|a|^2)^{(m+1)/2}$, while the Haar density contributes
$(1+|a|^2)^{-n/2}$. Their product is
\[
 (1+|a|^2)^{(m+1-n)/2}=(1+|a|^2)^{-\beta}.
\]
This proves the formula with exactly the stated constant.
The integrands are nonnegative and measurable in bundle charts, so Tonelli
justifies every rearrangement. The resulting integral is finite because the
section functional is finite for bounded $K$.
\end{proof}

\begin{remark}[The subcritical weight]
When $m\le n-2$, the weight decreases strictly as the slope length $|a|$
increases. It is this extra preference for small slopes that will exclude
noncentered fibers. At $m=n-1$ the two Jacobian powers cancel exactly, the
weight is constant, and the strictness proof below no longer applies.
\end{remark}

\subsection{From section moments to fiber constraints}

For $y\in u^\perp$ define the vertical fiber
$J_y(K)=\{s\in\R:y+su\in K\}$.
For $E\in\GG(u^\perp,m)$ and
$\mathbf y=(y_1,\ldots,y_{m+1})\in E^{m+1}$, set
\[
 R_{\mathbf y}(K)=J_{y_1}(K)\times\cdots\times J_{y_{m+1}}(K)
 \subset\R^{m+1},
 \qquad
 L_{\mathbf y}a=(\langle a,y_i\rangle)_{i=1}^{m+1}.
\]
For almost all $\mathbf y$, $L_{\mathbf y}:E\to\R^{m+1}$ has rank $m$.
Let
\[
 D_{\mathbf y}
 =\sqrt{\det\bigl(\sum_{i=1}^{m+1}y_i\otimes y_i\bigr)}
 \quad\text{(determinant on $E$)},
\]
let $\theta_{\mathbf y}\in\Sph^m$ span $\ker L_{\mathbf y}^*$,
so $L_{\mathbf y}(E)=\theta_{\mathbf y}^\perp$, and denote by
$L_{\mathbf y}^{-1}:\theta_{\mathbf y}^\perp\to E$ its inverse.

\begin{corollary}[Exact weighted fiber formula]\label{cor:fiber}
For any compact Borel set $K$,
\begin{equation}\label{eq:weighted-fiber}
 \begin{aligned}
 \Ph(K)=\frac{A_{n,m}}{m+1}
 \int_{\GG(u^\perp,m)}\int_{E^{m+1}}
 \frac{1}{D_{\mathbf y}}
 \int_{R_{\mathbf y}(K)\cap\theta_{\mathbf y}^\perp}
 \frac{\dd\HH^m(s)}
 {(1+|L_{\mathbf y}^{-1}s|^2)^\beta}
 \,\dd\mathbf y\,\dd\nu_{u^\perp,m}(E).
 \end{aligned}
\end{equation}
\end{corollary}
\begin{proof}
For fixed $(E,a)$, the power of a nonnegative integral is a product integral:
\[
 Q_K(E,a)^{m+1}
 =\int_{E^{m+1}}\prod_{i=1}^{m+1}
   \one_{J_{y_i}(K)}(\langle a,y_i\rangle)\,d\mathbf y.
\]
This is just repeated Tonelli; no symmetry or convexity of the fibers is
required. Interchanging the $a$ and $\mathbf y$ integrals gives the inner
quantity
\[
 \int_E w_\beta(a)\one_{R_{\mathbf y}(K)}(L_{\mathbf y}a)\,da.
\]
Suppose the tuple has rank $m$. Choose orthonormal coordinates on $E$ and
write $L=L_{\mathbf y}$ as an $(m+1)\times m$ matrix. Then
\[
 L^TL=\sum_{i=1}^{m+1}y_i\otimes y_i,
 \qquad J_mL=\sqrt{\det(L^TL)}=D_{\mathbf y}>0.
\]
The image is the $m$-plane $\theta_{\mathbf y}^{\perp}$, and
$d\HH^m(s)=D_{\mathbf y}\,da$ under $s=La$. Hence the last integral equals
\[
 \frac1{D_{\mathbf y}}
 \int_{R_{\mathbf y}(K)\cap\theta_{\mathbf y}^{\perp}}
          (1+|L_{\mathbf y}^{-1}s|^2)^{-\beta}\,d\HH^m(s).
\]
The deficient-rank tuples are the zero set of all maximal minors and have
Lebesgue measure zero. Applying the change of variables on full-rank tuples
proves the assertion. There is no need to estimate $D_{\mathbf y}^{-1}$
separately near rank deficiency: the equality is an identity of nonnegative
integrals, whose finiteness was already established by the graph formula.
\end{proof}

For estimates we will normally keep the integral in the $a$ variable rather
than use the factor $D_{\mathbf y}^{-1}$. In that form the conditions are
ordinary strips $\langle a,y_i\rangle\in J_{y_i}(K)$ and the weight is the
fixed smooth function $w_\beta$. The fiber formula and the strip formula
are exactly equivalent; the latter avoids an unnecessary apparent singularity.

The sign of $\theta_{\mathbf y}$ does not affect the integrand, which
can equivalently be written on $L_{\mathbf y}(E)$. On the rank-deficient
null set define the displayed integrand to be zero. This specifies a
measurable representative without choosing a sign for the kernel vector.

\begin{remark}[The classical endpoint]
For $m=n-1$, $\beta=0$ and $E=u^\perp$ is the sole element of
$\GG(u^\perp,n-1)$. Since
$A_{n,n-1}=\Gamma(n/2)/\pi^{n/2}=2/\sigma_n$,
\eqref{eq:weighted-fiber} reads
\[
 \Ph(K)=\frac{2}{n\sigma_n}
 \int_{(u^\perp)^n}D_{\mathbf y}^{-1}
 \HH^{n-1}(R_{\mathbf y}(K)\cap\theta_{\mathbf y}^\perp)
 \,\dd\mathbf y.
\]
Since $V(IK)=\sigma_n\Ph(K)$ at this endpoint, this is consistent
with the classical intersection-body volume representation in
\cite{MSY2025}.
\end{remark}

\begin{remark}[Normalization on balls]
For $K=RB_2^n$, the condition in $Q_K$ is
$|y|^2+\langle a,y\rangle^2\le R^2$, so
$Q_K(E,a)=\omega_mR^m(1+|a|^2)^{-1/2}$.
Substitution in \eqref{eq:weighted-graph} gives
$\Ph(RB_2^n)=\omega_m^{m+1}R^{m(m+1)}/(m+1)$,
using $A_{n,m}\int_{\R^m}(1+|a|^2)^{-n/2}\,da=1$.
\end{remark}

\section{Continuous Steiner symmetrization and weighted monotonicity}
\label{sec:weighted-monotonicity}

Throughout this section let $\beta=(n-m-1)/2\ge 0$, and let
$w_\beta(a)=(1+|a|^2)^{-\beta}$.  The pointwise formula
\eqref{eq:weighted-graph} will be used without any limiting singular
integrals.  All Grassmannians carry Haar probability measures.

\subsection{Finite-interval flows and the geometric input}
For a fixed $u$, write $J_y(K)=\{s:y+su\in K\}$ and
$D_u=P_{u^\perp}K$. The phrase \emph{$u$-multi-graphical} means that
there are disjoint open sets $\Omega_j\subset u^\perp$, $j\ge1$,
with $|D_u\setminus\bigcup_j\Omega_j|_{n-1}=0$, such that on
$\Omega_j$ the fiber is a union of exactly $j$ disjoint nondegenerate
closed intervals with ordered continuous endpoints. Put
$\Omega_\infty=\bigcup_j\Omega_j$.

For a finite interval union $J=\bigcup_i[c_i-\ell_i,c_i+\ell_i]$,
move each center to $(1-t)c_i$ until intervals first touch at time
$\tau$, merge touching components, and continue by
$S^tJ=S^{(t-\tau)/(1-\tau)}(S^\tau J)$. This prescription has finitely
many mergers. The compact-set construction is
\[
 K_t=S_u^tK=\overline{\{y+su:y\in\Omega_\infty,\ s\in S^tJ_y(K)\}}.
\]
We use the following published geometric facts, with numbering from
the version of record of~\cite{MSY2025}.

\begin{proposition}[Geometric input from MSY]\label{prop:MSY-input}
For the construction just described, the compact set $K_t$ is independent
of the graphical representation and
\[
 J_y(K_t)=S^tJ_y(K)\quad(y\in\Omega_\infty,\ 0\le t\le1).
\]
Volume and inclusion between compact multi-graphical sets are preserved.
For a Lipschitz star body, almost every $u$ has the stated graphical representation; along every such
direction every $K_t$ is a Lipschitz star body, $K_0=K$, $K_1=S_uK$,
and for some $M_0<\infty$,
\begin{equation}\label{eq:MSY-radial}
 (1+M_0t)^{-1}\le\frac{\rho_{K_t}(\theta)}{\rho_K(\theta)}
 \le1+M_0t\qquad(\theta\in\Sph^{n-1},\ 0\le t\le1).
\end{equation}
There is also $\delta>0$ with $B_2^n(\delta)\subset\operatorname{int}K$
such that the fibers over $|y|<\delta$ are single nondegenerate intervals
with continuous endpoints.
\end{proposition}
\begin{proof}[Published inputs]
Representation independence and the simultaneous fiber identity are
\cite[Proposition 4.8 and Corollary 4.10]{MSY2025}; volume and inclusion
are Corollaries 4.12 and 4.14. The remaining statements are
Proposition 5.6, Theorem 5.7, Proposition 5.8, Corollary 5.9, and
Lemma 5.11 there. All apply to Lipschitz star bodies, without convexity
or any fixed-point hypothesis. No strictness or equality classification
from that paper is used here.
\end{proof}

For arbitrary compact multi-graphical sets, endpoint equalities may hold only
up to volume-null sets; exact $K_0=K$ and $K_1=S_uK$ above are asserted for
Lipschitz star bodies. Also, full-measure families of directions need not
contain every direction. The argument first proves symmetry for those good
directions and only later extends it to all directions by continuity.

A finite union of intervals evolves in a simple way between collisions:
lengths stay fixed and all centers are multiplied by the same scalar.
After a collision the number of components drops. Consequently a fixed
finite collection of fibers requires only finitely many collision-free
time intervals. The number of components need not be uniformly bounded
over every base point; such a global bound is not used.

\begin{lemma}[A common star-shapedness ball and radial control]
\label{lem:radial-core}
A Lipschitz star body is star-shaped with respect to every point of
some ball $\delta B_2^n$, $\delta>0$. If two star bodies $K,L$ share
such a ball and $d_H(K,L)\le d$, then
\[
 (1+d/\delta)^{-1}\le\frac{\rho_L}{\rho_K}\le1+d/\delta.
\]
\end{lemma}
\begin{proof}
The gauge $p_K(x)=|x|/\rho_K(x/|x|)$, with $p_K(0)=0$, is globally
Lipschitz. To see this, write $h=1/\rho_K$ and $x=ru$, $y=sv$,
where $r\ge s>0$. Then
\[
 |rh(u)-sh(v)|
 \le\|h\|_\infty|r-s|+\operatorname{Lip}(h)s|u-v|
 \le(\|h\|_\infty+2\operatorname{Lip}(h))|x-y|.
\]
The case $y=0$ is immediate. Let $L_K$ be a Lipschitz constant and
choose $L_K\delta<1$. For $x\in K$, $z\in\delta B_2^n$ and
$0\le t<1$,
\[
 p_K(tx+(1-t)z)\le t p_K(x)+(1-t)L_K|z|<1.
\]
Thus every segment from $z$ to $x$ lies in $K$.

If $K$ has this property, then
$K+dB_2^n\subset(1+d/\delta)K$.
For $d>0$, put $\epsilon=d/\delta$ and observe that
\[
 \frac{x+z}{1+\epsilon}
 =\frac{x}{1+\epsilon}+\frac{\epsilon}{1+\epsilon}\frac z\epsilon
 \in K\qquad(x\in K,\ |z|\le d).
\]
The case $d=0$ is immediate. Apply this inclusion to both Hausdorff
inclusions $L\subset K+dB_2^n$ and $K\subset L+dB_2^n$ and then
compare radial functions.
\end{proof}

Along the flow, preservation of a common star-shapedness ball and the
bound $d_H(K_t,K)\le Rt$ for $K\subset RB_2^n$ are
\cite[Proposition 4.16 and Lemma 4.17]{MSY2025}, together with its
Corollary 5.9. Lemma~\ref{lem:radial-core} therefore also explains the
radial estimate \eqref{eq:MSY-radial}. Merely containing an inner ball
would not suffice for this comparison.

\subsection{The weighted inequality}

For arbitrary nonorthogonal strip constraints, translating their center
vector toward zero is not a separate coordinate motion. We therefore do
not infer monotonicity just from radial decrease of $w_\beta$.
The proof first treats Gaussian weights, for which joint log-concavity is
available, and then averages those inequalities with positive coefficients.
A positive mixture of log-concave functions need not be log-concave; only
their common one-dimensional monotonicity is used.

\begin{lemma}[Even log-concavity of Gaussian box integrals]
\label{lem:gaussian}
Let $E$ be an $m$-dimensional Euclidean space, let
$y_1,\ldots,y_{m+1}\in E$ span $E$, and let $\ell_i>0$.  For
$\gamma\ge0$ and $c=(c_i)\in\mathbb R^{m+1}$, put
\[
 H_\gamma(c)=\int_E e^{-\gamma|a|^2}
 \prod_{i=1}^{m+1}\one_{[-\ell_i,\ell_i]}
 (\langle a,y_i\rangle-c_i)\,\dd a.
\]
Then $H_\gamma$ is finite, even and log-concave.  Hence
\[
 H_\gamma(sc)\ge H_\gamma(tc)
 \qquad (0\le s\le t,\ c\in\mathbb R^{m+1}).
\]
\end{lemma}
\begin{proof}
The spanning assumption makes every feasible intersection of strips
bounded, so $H_\gamma(c)<\infty$ even for $\gamma=0$. For $\gamma>0$
the function
\[
 (a,c)\longmapsto e^{-\gamma|a|^2}
 \one_{\{|\langle a,y_i\rangle-c_i|\le\ell_i\ \forall i\}}
\]
is log-concave: its constraint set is convex. It is also integrable
on the whole product space, since its integral there equals
$\prod_i(2\ell_i)\int_Ee^{-\gamma|a|^2}\dd a<\infty$.
The marginal theorem of Pr\'ekopa~\cite[Theorem 6]{Prekopa1973} says that
integrating a nonnegative integrable log-concave function over some of its
variables leaves a log-concave function of the remaining variables. All
its hypotheses have just been checked, so $H_\gamma$ is log-concave. For $\gamma\downarrow0$, monotone
convergence gives $H_\gamma(c)\uparrow H_0(c)<\infty$; passing to the
limit in the log-concavity inequality proves it also for $H_0$.
Changing $a$ to $-a$ gives evenness. For $0\le s<t$, express
$sc=\frac{1+s/t}{2}tc+\frac{1-s/t}{2}(-tc)$ and apply evenness and
log-concavity. The endpoint cases $s=t$ and $t=0$ are identities.
\end{proof}

\begin{proposition}[Weighted rectangle monotonicity]
\label{prop:rect}
Let $E$, $y_i$, $\ell_i$ be as in Lemma~\ref{lem:gaussian},
$\beta\ge0$, and $R=\prod_{i=1}^{m+1}[c_i-\ell_i,c_i+\ell_i]$.
For $0\le t\le1$ define
\[
 \mathcal J_R(t)=\int_E w_\beta(a)
 \one_{R-tc}\bigl(L_{\mathbf y}a\bigr)\,\dd a,
 \qquad L_{\mathbf y}a=(\langle a,y_i\rangle)_{i=1}^{m+1}.
\]
Then $\mathcal J_R(t)$ is non-decreasing on $[0,1]$.
\end{proposition}
\begin{proof}
For $\beta=0$ this is exactly Lemma~\ref{lem:gaussian} with
$\gamma=0$, since $R-tc$ has center $(1-t)c$.  For $\beta>0$ use
\[
 w_\beta(a)=\frac1{\Gamma(\beta)}
   \int_0^\infty\gamma^{\beta-1}e^{-\gamma}
             e^{-\gamma|a|^2}\,\dd\gamma.
\]
For $0\le t_1\le t_2\le1$, Lemma~\ref{lem:gaussian} gives
\[
 H_\gamma((1-t_2)c)\ge H_\gamma((1-t_1)c).
\]
Integrate this inequality against the nonnegative density
$\gamma^{\beta-1}e^{-\gamma}/\Gamma(\beta)$ and use Tonelli.
The mixture identity itself follows from the Gamma integral after the
substitution $s=\gamma(1+|a|^2)$. This proves the assertion.
\end{proof}

\subsection{A single exceptional set for the whole flow}

The graphical description is valid outside a base-space null set.
The following incidence identity transfers this set to the parameter
bundle and supplies a single exceptional set for the entire flow.

\begin{lemma}[Null-set transfer to the parameter bundle]
\label{lem:common-good-set}
Fix a compact $u$-multi-graphical $K$. In the bundle
\[
 \mathcal P_u=\{(E,y_1,\ldots,y_{m+1}):E\in\GG(u^\perp,m),\ y_i\in E\},
 \qquad \dd\mu=\dd\nu(E)\,\dd y_1\cdots\dd y_{m+1},
\]
there is a single $\mu$-conull measurable set, independent of $t$, on
which all $y_i$ are nonzero, span $E$, and either every $y_i$ belongs
to $\Omega_\infty$, or at least one $y_i$ lies outside $D_u$.
In the first case the exact finite-interval identity holds for all
times; in the second case the product of the fiber indicators is zero
for all times.
\end{lemma}
\begin{proof}
Put $d=n-1$. Polar coordinates in each $E$, followed by rotational
invariance of the angular incidence probability, give for $g\ge0$
\begin{equation}\label{eq:polar-incidence}
 \int_{\GG(d,m)}\int_E g(y)\,\dd y\,\dd\nu(E)
 =\frac{\sigma_m}{\sigma_d}\int_{\R^d}g(y)|y|^{m-d}\,\dd y.
\end{equation}
In detail, the left side is $\sigma_m\int_0^\infty r^{m-1}$ times
the average of $g(r\theta)$ over incident $(E,\theta)$; the marginal
of $\theta$ is uniform on $\Sph^{d-1}$, which is the right side in
polar coordinates. The density $|y|^{m-d}$ is locally integrable
because its radial integral near zero is $\int_0^1r^{m-1}\dd r<\infty$.
Apply this identity to the indicator of the null set
$N=D_u\setminus\Omega_\infty$, truncated in balls of integer radius.
For almost every $E$, $N\cap E$ is null in $E$. Fiberwise Fubini then
excludes tuples having any coordinate in $N$. Tuples with a zero
coordinate or deficient rank are also null, by ordinary Lebesgue
measure and the vanishing of nontrivial minors. These exclusions do
not involve $t$. The simultaneous identity on $\Omega_\infty$ is
Proposition~\ref{prop:MSY-input}. Finally $P_{u^\perp}K_t\subset D_u$
for all $t$, because $D_u$ is closed and the construction uses only
base points in $D_u$. A coordinate outside $D_u$ therefore makes
the product vanish for every $t$.
\end{proof}

\begin{theorem}[Monotonicity under continuous Steiner symmetrization]
\label{thm:monotonicity}
Let $K$ be a compact $u$-multi-graphical set in the sense of
\cite[Definition 4.2]{MSY2025}, and let $K_t=S_u^tK$ denote the
finite-interval collision-flow symmetrization of
\cite[Section 4.2]{MSY2025}.  Then for $1\le m\le n-1$,
\[
 [0,1]\ni t\longmapsto\Ph(K_t)
 \quad\text{is non-decreasing}.
\]
\end{theorem}
\begin{proof}
The graph formula and Tonelli give, for every $t$,
\begin{equation}\label{eq:rect-kernel}
 \Ph(K_t)=\frac{A_{n,m}}{m+1}
 \int_{\GG(u^\perp,m)}\int_{E^{m+1}}
       \mathcal J_{E,\mathbf y}(t)\,\dd\mathbf y\,\dd\nu(E),
\end{equation}
where
\[
 \mathcal J_{E,\mathbf y}(t)=\int_E w_\beta(a)
       \prod_{i=1}^{m+1}\one_{J_{y_i}(K_t)}(\langle a,y_i\rangle)\,\dd a.
\]
Work on the single conull set of Lemma~\ref{lem:common-good-set}.
If a coordinate lies outside the projection, $\mathcal J$ is identically
zero. Otherwise all the involved fibers are exact finite-interval
flows. Their combined collision times are finite in number.
On a collision-free interval starting at $\tau$, their product is
a finite union of rectangles $R_b(t)$ with disjoint interiors, where
\[
 R_b(t)=R_b(\tau)-\frac{t-\tau}{1-\tau}c_b(\tau)
\]
and $c_b(\tau)$ is the center of $R_b(\tau)$.
Indeed, a component created at an earlier time $s\le\tau$ has center
$C(t)=(1-t)C(s)/(1-s)$ until its next collision. Therefore
$C(t)=(1-t)C(\tau)/(1-\tau)$ on the current interval, regardless of
its creation time. All the selected fibers consequently use the same
contraction parameter, even when their collision times differ.
Each rectangle contributes a non-decreasing integral by
Proposition~\ref{prop:rect}.

The finite rectangle contributions add without overlap in the integral:
a common facet is given by $\langle a,y_i\rangle=s$, a null hyperplane
in $E$ because $y_i\ne0$. At a collision only the decomposition changes;
the interval unions have continuous endpoints. Outside the finitely
many corresponding boundary hyperplanes the product indicators converge
pointwise as $t$ approaches the collision. All endpoints remain uniformly
bounded over $[0,1]$. Since the $y_i$ span $E$, their strips confine $a$
to a fixed ball for this tuple. Thus bounded convergence proves continuity
of $\mathcal J$ at each collision, and also at $t=1$. Concatenating the
finitely many monotone pieces proves that $\mathcal J(t)$ is
non-decreasing on $[0,1]$ on the same conull parameter set.
Integrating the pointwise comparison proves the theorem.
\end{proof}

\begin{remark}
The collision-flow construction and its measurable graphical-fiber
identity are used from~\cite{MSY2025}; the new part here is the
Gaussian-mixture argument for the weight $w_\beta$.  No strictness
follows from the preceding log-concavity assertion by itself.
\end{remark}

\section{Strict growth and fixed-point rigidity}
\label{sec:weighted-strict}

Assume $1\le m\le n-2$, so that $\beta=(n-m-1)/2>0$.
We strengthen monotonicity to a uniform positive linear increment
whenever the initial body is not Steiner symmetric. The quadratic
estimate of Proposition~\ref{prop:quadratic-stationarity} will then
force every Lipschitz fixed point to be a ball.

\subsection{Strictness for a product of intervals}

Nondecreasing functions can have zero derivative even when they are not
constant. Therefore monotonicity alone cannot contradict stationarity.
We must exhibit a particular family of strip configurations with a strictly
positive initial rate, make that rate uniform on a positive-measure set,
and then integrate finite increments. The one-variable calculation below
is the sole source of strictness.

\begin{lemma}[Strict motion of one interval]\label{lem:1d-strict}
Let $\beta>0$, $B\ge0$, $r,\ell>0$, and $c\in\R$. Set
$W_B(z)=(1+B+z^2)^{-\beta}$ and
\[
 F(t)=\int_{((1-t)c-\ell)/r}^{((1-t)c+\ell)/r}W_B(z)\,\dd z.
\]
For $0\le t<1$,
\begin{equation}\label{eq:1d-deriv}
 F'(t)=\frac c r\left[
 W_B\left(\frac{(1-t)c-\ell}{r}\right)
 -W_B\left(\frac{(1-t)c+\ell}{r}\right)\right]\ge0,
\end{equation}
and equality holds if and only if $c=0$.
\end{lemma}
\begin{proof}
Set $L(t)=((1-t)c-\ell)/r$ and $U(t)=((1-t)c+\ell)/r$.
Both derivatives are $-c/r$, so Leibniz' rule gives
\[
 F'(t)=W_B(U(t))\frac{-c}{r}-W_B(L(t))\frac{-c}{r},
\]
which is the stated identity. To check the sign without a case ambiguity,
observe that
\[
 ((1-t)c+\ell)^2-((1-t)c-\ell)^2=4(1-t)c\ell.
\]
For $c>0$ and $t<1$, $|L(t)|<|U(t)|$, and the bracket is strictly
positive. For $c<0$ the bracket is strictly negative and is multiplied by
a negative $c/r$. In either nonzero case the product is strictly positive.
For $c=0$ the interval does not move and $F'(t)=0$.
Although $W_B'(0)=0$, $W_B$ is strictly decreasing as a function of
$|z|$ on $[0,\infty)$; the comparison of two distinct absolute values is
therefore strict. The endpoint $t=1$ is deliberately excluded.
\end{proof}

\begin{lemma}[A fixed-cube formula and parameter continuity]
\label{lem:fixed-cube}
Let $B\in\mathrm{GL}(m,\R)$, $c\in\R^m$, and
$D=\operatorname{diag}(\ell_1,\ldots,\ell_m)$ with $\ell_i>0$.
Define
\[
 I(B,c,D;t)=\int_{\R^m}w_\beta(a)
 \one_{(1-t)c+D[-1,1]^m}(Ba)\,\dd a.
\]
Then
\begin{equation}\label{eq:fixed-cube}
 I(B,c,D;t)=\frac{\det D}{|\det B|}
 \int_{[-1,1]^m}
 w_\beta\bigl(B^{-1}((1-t)c+Dv)\bigr)\,\dd v.
\end{equation}
Both $I$ and $\partial_t I$ are jointly continuous in $(B,c,D,t)$.
More explicitly, with $a(v,t)=B^{-1}((1-t)c+Dv)$,
\begin{equation}\label{eq:fixed-cube-derivative}
 \partial_t I
 =\frac{2\beta\det D}{|\det B|}
 \int_{[-1,1]^m}
 \frac{\langle a(v,t),B^{-1}c\rangle}
 {(1+|a(v,t)|^2)^{\beta+1}}\,\dd v.
\end{equation}
\end{lemma}
\begin{proof}
Use the affine substitution $Ba=(1-t)c+Dv$. On a compact
neighborhood in the displayed parameter domain, $B^{-1}$,
$|\det B|^{-1}$, $c$, and $D$ are uniformly bounded. The arguments
of $w_\beta$ remain bounded for bounded $t$, and the integrands in
\eqref{eq:fixed-cube} and \eqref{eq:fixed-cube-derivative} are
continuous and uniformly bounded on the fixed cube. Differentiation
in $t$ and joint continuity follow by dominated convergence.
Notice that no derivatives of the parameters $c$ or $D$ with
respect to other variables are needed.
\end{proof}

\begin{remark}[Nonzero centers alone do not imply strictness]
The numerator in \eqref{eq:fixed-cube-derivative} need not be
pointwise positive. In the overdetermined $(m+1)$-strip formula,
a nonzero center can even belong to an inactive constraint.
For example, in $\R^2$ take $y_1=e_1$, $y_2=e_2$,
$y_3=e_1+e_2$, half-lengths $(1,1,10)$ and centers $(0,0,1)$.
For $0\le t\le1$ the third strip is redundant on $[-1,1]^2$,
so the weighted integral is constant although the center vector
is nonzero. This is not a counterexample to Lemma~\ref{lem:strict}:
the proof below deliberately makes a \emph{noncentered active}
constraint a coordinate interval.
\end{remark}

\subsection{Uniform localization in the parameter bundle}

Use the bundle $\mathcal P_u$ and measure $\mu$ from
Lemma~\ref{lem:common-good-set}, and write
$\mathcal J_P(t)=\mathcal J_{E,\mathbf y}(t)$ for the integrand in
\eqref{eq:rect-kernel}.

\begin{lemma}[Adjoining the central graphical neighborhood]
\label{lem:central-compatible}
Let $K$ be a Lipschitz star body and a $u$-multi-graphical set.
Its graphical representation can be chosen so that
$B_{u^\perp}(\delta)\subset\Omega_1$ for some $\delta>0$ with
$B_2^n(\delta)\subset\operatorname{int}K$, without changing any $K_t$.
Consequently the exact finite-interval identity holds at every point
of this central neighborhood, including $y=0$, for every $t\in[0,1]$.
\end{lemma}
\begin{proof}
By Proposition~\ref{prop:MSY-input}, after decreasing $\delta$ if
necessary, $U_0=B_{u^\perp}(\delta)$ has single nondegenerate fibers
$[f_0(y),g_0(y)]$ with continuous endpoints. Starting with any graphical
representation, put $\Omega'_1=\Omega_1\cup U_0$ and
$\Omega'_j=\Omega_j$ for $j\ge2$. The sets are disjoint: a fiber
cannot be both one interval and two or more disjoint intervals separated
by positive gaps. On $U_0\cap\Omega_1$ the ordered endpoints agree,
so the endpoint functions glue continuously. The uncovered part of the
projection can only decrease. This is therefore another valid graphical
representation.

For completeness, independence for this particular enlargement follows
directly, without requiring any inference from a full-measure set about
a prescribed point. The old unclosed flow is contained in the new one.
For $y\in U_0$, choose $y_j\to y$ with
$y_j\in U_0\cap\Omega_\infty$, possible because the latter has full
measure in the open set $U_0$. Such $y_j$ must belong to $\Omega_1$.
For each $t$, the single-interval endpoints are
$(1-t)c_0(y)\pm\ell_0(y)$, where
$c_0=(f_0+g_0)/2$ and $\ell_0=(g_0-f_0)/2$.
They are continuous in $y$. Every point of the added fiber is thus
a limit of points of the old unclosed flow, by keeping its relative
position inside the interval fixed. The two closures agree for every
$t$. Corollary 4.10 of~\cite{MSY2025}, applied to the enlarged valid
representation, gives the exact fiber identity there. Relabel this
representation as $\{\Omega_j\}$.
\end{proof}

\begin{lemma}[Uniform positive growth unless the fibers are centered]
\label{lem:strict}
Let $K$ be a Lipschitz star body and let $u$ be a direction for which
$K$ is $u$-multi-graphical. Let $K_t=S_u^tK$ be the finite-interval
continuous Steiner symmetrization used in Section~\ref{sec:weighted-monotonicity}.
If $K$ does not agree almost everywhere with $S_uK$, then there
exist $\varepsilon>0$ and $\kappa>0$ such that
\begin{equation}\label{eq:strict-finite}
 \Ph(K_t)-\Ph(K)\ge\kappa t
 \qquad (0<t\le\varepsilon).
\end{equation}
In particular, the lower right Dini derivative is strictly positive.
\end{lemma}
\begin{proof}
\textbf{1. Fixed graphical neighborhoods.}
Use Lemma~\ref{lem:central-compatible} to arrange
$U_0=B_{u^\perp}(\delta)\subset\Omega_1$ and
$B_2^n(\delta)\subset\operatorname{int}K$.
The geometric input preserves this ball for all $t$.
If every good fiber consisted of a single interval centered at zero,
then Fubini would give $|K\mathbin{\triangle}S_uK|=0$.
A union of two or more disjoint nondegenerate intervals cannot have
all component centers equal to zero. The hypothesis and Fubini therefore
supply $y_1^0\ne0$ in some $\Omega_j$ with a component of nonzero center. Shrink a neighborhood
of $y_1^0$ with compact closure in $\Omega_j$ so that its finite list
of half-lengths and gaps has positive lower bounds and the selected
center remains nonzero.

Choose $E_0\in\GG(u^\perp,m)$ containing $y_1^0$ and complete it
to a nonzero orthogonal basis $y_1^0,\ldots,y_m^0$, taking
$|y_i^0|<\delta/4$ for $i\ge2$. Put $y_{m+1}^0=0$ and
$P_0=(E_0,\mathbf y^0)$. When $m=1$ there are no extra basis vectors.
Every central fiber used here belongs to the enlarged $\Omega_1$;
no assertion that a prescribed vector lies in an arbitrary conull set
is needed.

\textbf{2. Strict derivative at the anchor.}
Set $r_i=|y_i^0|$ and write $a=\sum_{i=1}^m z_i y_i^0/r_i$.
The last constraint is automatic, since $0\in J_0(K_t)$. Hence
\begin{equation}\label{eq:orthogonal-anchor}
 \mathcal J_{P_0}(t)=\int_{\R^m}(1+|z|^2)^{-\beta}
       \prod_{i=1}^m\one_{J_{y_i^0}(K_t)}(r_iz_i)\,\dd z.
\end{equation}
Before the first collision, decompose this integral into finitely many
product boxes, one for each component of $J_{y_1^0}(K)$; the other
fibers are single intervals. For one box write its coordinate endpoints as
\[
 L_i(t)=\frac{(1-t)c_i-\ell_i}{r_i},\qquad
 U_i(t)=\frac{(1-t)c_i+\ell_i}{r_i}.
\]
Its integral $G(t)$ has derivative
\[
 G'(t)=\sum_{i=1}^{m}\frac{c_i}{r_i}
 \int_{\prod_{h\ne i}[L_h(t),U_h(t)]}
 \left[
 (1+\sum_{h\ne i}z_h^2+L_i(t)^2)^{-\beta}
 -(1+\sum_{h\ne i}z_h^2+U_i(t)^2)^{-\beta}
 \right]d\widehat z_i.
\]
Here $d\widehat z_i$ means integration in every coordinate except $i$.
This is obtained by differentiating one pair of moving endpoints at a
time and adding the contributions; it does not differentiate the full
Grassmannian integral. Every interval is bounded, so ordinary Leibniz
rules and Fubini apply. Lemma~\ref{lem:1d-strict}, with $B=\sum_{h\ne i}z_h^2$, shows that
each contribution is nonnegative. For the selected first component,
$c_1\ne0$, so its first-coordinate contribution is positive for every
choice of the other coordinates. Those coordinates vary over intervals
of positive length. When $m=1$, there are no other coordinates: integration over $\R^0$ leaves the displayed integrand unchanged, since the empty product domain has measure $1$. Therefore
\begin{equation}\label{eq:anchor-positive}
 d_0:=\partial_t\mathcal J_{P_0}(0)>0.
\end{equation}
All these integrals and derivatives are finite because their domains
are bounded.

\textbf{3. One parameter neighborhood and one time interval.}
Choose a smooth orthonormal frame of $E$ near $E_0$ and use it to
trivialize the bundle. Restrict the first point to the graphical
neighborhood just chosen, the next $m-1$ points to small neighborhoods
inside $U_0$, and keep the first $m$ vectors linearly independent.
Let $B_P$ be the matrix of
$a\mapsto(\langle a,y_1\rangle,\ldots,\langle a,y_m\rangle)$.
On a sufficiently small relatively compact neighborhood,
$\|B_P^{-1}\|\le C_B$. Take $R\ge1$ with $K\subset B_2^n(R)$.
All initial centers and endpoints have absolute value at most $R$.
If the first fiber has more than one component, let $g_*>0$ be a
common lower bound for its initial gaps. Before collision the gap
between adjacent intervals equals its initial gap minus
$t(c_{i+1}-c_i)$, and $|c_{i+1}-c_i|\le2R$. Thus
\[
 \varepsilon_0=\min\{1/2,g_*/(4R)\}>0
\]
is a common collision-free time; when there are no gaps take
$\varepsilon_0=1/2$. The other chosen fibers each have one component.
Throughout this interval their endpoints remain in $[-R,R]$.
The first $m$ strips therefore imply
\[
 |a|\le C_B\sqrt m\,R\le M
\]
with one fixed $M\ge1$ for all nearby parameters and times.

Now additionally impose
$|y_{m+1}|<\min\{\delta/4,\delta/(4M)\}$.
On the domain determined by the first $m$ constraints,
$|\langle a,y_{m+1}\rangle|<\delta/4$.
Every point $y_{m+1}+su$ with $|s|\le\delta/4$ then lies in
$B_2^n(\delta)\subset K_t$. The last constraint is consequently
redundant on this whole open parameter neighborhood and throughout
$[0,\varepsilon_0]$.

\textbf{4. Uniform strictness on positive measure.}
On this neighborhood $\mathcal J_P(t)$ is a fixed finite sum of
\[
 I(B_P,c_b(P),D_b(P);t),
\]
where $b$ indexes the components of the first fiber. Its centers and
half-length matrices are continuous in $P$. Lemma~\ref{lem:fixed-cube}
therefore makes $\partial_t\mathcal J_P(t)$ jointly continuous in
$(P,t)$ on the common collision-free interval. By
\eqref{eq:anchor-positive}, there are a relatively compact open
neighborhood $\mathcal U$ of $P_0$ and $\varepsilon\in(0,\varepsilon_0]$
such that
\begin{equation}\label{eq:uniform-positive}
 \partial_t\mathcal J_P(t)\ge d_0/2
 \qquad(P\in\mathcal U,\ 0\le t\le\varepsilon).
\end{equation}
In these coordinates Haar measure has positive smooth density and
fiber measure is ordinary Lebesgue measure. Thus
$0<\mu(\mathcal U)<\infty$. Time integration gives
$\mathcal J_P(t)-\mathcal J_P(0)\ge(d_0/2)t$ on $\mathcal U$.

\textbf{5. Integration of finite increments.}
On the common conull parameter set from
Lemma~\ref{lem:common-good-set}, the proof of
Theorem~\ref{thm:monotonicity} gives
$\mathcal J_P(t)-\mathcal J_P(0)\ge0$ for every $t$.
The functionals in \eqref{eq:rect-kernel} are finite, and $K_0=K$
by Proposition~\ref{prop:MSY-input}. Subtraction and integration are
therefore legitimate and yield
\[
 \Ph(K_t)-\Ph(K)
 \ge\frac{A_{n,m}}{m+1}\frac{d_0}{2}\mu(\mathcal U)t.
\]
Taking $\kappa=A_{n,m}d_0\mu(\mathcal U)/(2(m+1))>0$ proves the
claim. No derivative or limit is moved through the outer integral.
\end{proof}

\begin{remark}[The set on which strictness is integrated]
The orthogonal anchor tuple and $y_{m+1}=0$ were only used to find one
strictly positive derivative. They were never integrated as a set of
positive measure. The actual integration takes place on an open set in
all variables, including nonorthogonal tuples and nonzero $y_{m+1}$.
Strictness survives there by joint continuity of the derivative, not by
an assertion that every nonorthogonal configuration is strictly monotone.

In a local trivialization, $\mathcal P_u$ has ordinary coordinates of
dimension $m(n-1-m)+m(m+1)=mn$. The chosen set $\mathcal U$ is open in
these coordinates. Haar probability has a positive smooth local density,
and the remaining factors are Lebesgue measures, explaining directly
why $\mu(\mathcal U)>0$.

\end{remark}

\subsection{The fixed-point classification}

We now apply the two increment estimates to the same volume-preserving
path. The contradiction concerns finite differences, not a formal
Euler--Lagrange calculation along the deformation.

\begin{corollary}[Fixed-point rigidity for Lipschitz star bodies]
\label{cor:lipschitz}
Let $n\ge3$, $1\le m\le n-2$, and let $K$ be a Lipschitz star body
satisfying $\mathbb I_mK=cK$ for some $c>0$. Then $K$ is an
origin-centered Euclidean ball.
\end{corollary}
\begin{proof}
By \cite[Theorem 5.7]{MSY2025}, $K$ is $u$-multi-graphical for a
full-measure set of directions $u$. For each such $u$, Proposition~\ref{prop:MSY-input} gives star bodies
$K_t=S_u^tK$ with $K_0=K$ and $|K_t|=|K|$.
From \eqref{eq:MSY-radial}, with $R=\max\rho_K$, both signs of the
radial difference satisfy
\[
 \|\rho_{K_t}-\rho_K\|_\infty\le M_0Rt=:Mt.
\]
Indeed the downward difference is at most
$R[1-(1+M_0t)^{-1}]\le RM_0t$, and the upward difference is at most
$RM_0t$.
Proposition~\ref{prop:quadratic-stationarity} gives
$|\Ph(K_t)-\Ph(K)|\le C_u t^2$ for small $t$.
If $K\ne S_uK$ up to null sets, Lemma~\ref{lem:strict} gives
$\Ph(K_t)-\Ph(K)\ge\kappa_u t$ with $\kappa_u>0$.
Choose $t>0$ smaller than both validity times and, in addition,
$t<\kappa_u/(2(C_u+1))$. The two inequalities would then imply
$\kappa_u\le C_ut<\kappa_u/2$, a contradiction. Thus they cannot
hold for any non-Steiner-symmetric fixed point.
Hence $K=S_uK$ almost everywhere for every direction in this
full-measure set. In particular, $K$ agrees almost everywhere with
its reflection $R_uK$ about $u^\perp$.

The identity
\[
 |K\mathbin{\triangle}R_uK|
 =\frac1n\int_{\Sph^{n-1}}
   |\rho_K(\theta)^n-\rho_{R_uK}(\theta)^n|\,\dd\HH^{n-1}(\theta)
\]
shows that the nonnegative continuous integrand has integral zero.
A continuous nonnegative function that is positive at one point is
positive on an open set of positive spherical measure, so that is
impossible. Hence the radial functions agree at every point and
$K=R_uK$ exactly.
Every nonempty spherical open set has positive measure, so the good
directions are dense. For an arbitrary $u$, choose good directions $u_j\to u$.
For each $\theta$, the reflection formula
$R_{u_j}\theta=\theta-2\langle\theta,u_j\rangle u_j$ converges to $R_u\theta$.
Continuity of $\rho$ therefore gives
$\rho(R_u\theta)=\lim_j\rho(R_{u_j}\theta)=\rho(\theta)$.
The reflection identity thus holds for every direction.
For any distinct $v,w\in\Sph^{n-1}$, the reflection with normal
$(v-w)/|v-w|$ sends $v$ to $w$. Therefore $\rho_K(v)=\rho_K(w)$,
and $K$ is a centered ball.
\end{proof}

\begin{theorem}[Fixed points of the bi-dual $m$-intersection operator]
\label{thm:global-fixed}
Let $n\ge3$ and $1\le m\le n-2$, and let $K$ be a star body in
$\R^n$. Then
\[
 \mathbb I_mK=cK\quad\text{for some }c>0
 \quad\Longleftrightarrow\quad
 K=rB_2^n\quad\text{for some }r>0.
\]
\end{theorem}
\begin{proof}
For $m=1$ the assertion follows from
Proposition~\ref{prop:m1}. Assume $2\le m\le n-2$.
Theorem~\ref{thm:smooth} implies that every fixed-point star body
has a smooth, and thus Lipschitz, radial function. Consequently
Corollary~\ref{cor:lipschitz} applies and forces the body to be an
origin-centered Euclidean ball.
Conversely, for $K=rB_2^n$ we have
$\vol_m(K\cap\xi)=\omega_mr^m$ for every $\xi$, so
\[
 \rho_{\mathbb I_mK}^{n-m}(u)
 =\frac{\sigma_m}{\sigma_n}\omega_m^m r^{m^2},
\]
independently of $u$. Therefore $\mathbb I_mK=cK$ with
\[
 c=\left(\frac{\sigma_m}{\sigma_n}\omega_m^m\right)^{1/(n-m)}
 r^{\frac{m^2}{n-m}-1}>0.
\]
\end{proof}

For a ball, the displayed value of $c$ is positive for every radius.
The classification therefore concerns shape rather than a distinguished
unnormalized radius. In particular, it makes no assertion that all balls
are exact fixed points of the unnormalized operator.

\section{A strict Lyapunov functional and periodic points}\label{sec:periodic}

The argument in this section does not use regularity or symmetrization.
We first establish an abstract norm inequality and then apply it to the
probability Radon averages. It is the precise equality condition, rather
than monotonicity alone, that excludes nontrivial periodic shapes.

\begin{lemma}[Strict monotonicity for a positive nonlinear power step]
\label{lem:nonlinear-power}
Let $(X,\mu)$ and $(Y,\nu)$ be probability spaces,
let $1<p,q<\infty$, and let $A$ be a positive linear operator
with an adjoint $A^*$ satisfying
$\int_Y(Ah)k\,d\nu=\int_Xh(A^*k)\,d\mu$
for the positive functions involved below.
Suppose $f>0$ and
\[
 g=A^*\big[(Af)^{q-1}\big]>0,\qquad
 f_+=g^{1/(p-1)}.
\]
Assume the indicated functions and their images have the finite
positive norms used below. Then
\begin{equation}\label{eq:nonlinear-power-increase}
 \frac{\|Af_+\|_q}{\|f_+\|_p}
 \ge \frac{\|Af\|_q}{\|f\|_p},
\end{equation}
with equality if and only if $f_+$ is a positive scalar multiple of
$f$ almost everywhere.
\end{lemma}
\begin{proof}
Set $p'=p/(p-1)$ and $q'=q/(q-1)$. All norms below are finite and
positive by hypothesis, so all divisions are legitimate.

\textbf{1. Compare the new numerator to a dual pairing.}
H\"older on $Y$ with exponents $q,q'$ gives
\[
 \int_Y(Af_+)(Af)^{q-1}\,d\nu
 \le\|Af_+\|_q\,\|(Af)^{q-1}\|_{q'}
 =\|Af_+\|_q\,\|Af\|_q^{q-1}.
\]
Adjointness turns the left side into $\int_Xf_+g\,d\mu$.
By the specific choice $f_+=g^{1/(p-1)}$,
\[
 \int_X f_+g\,d\mu=\int_X g^{p'}\,d\mu=\|g\|_{p'}^{p'},
 \qquad \|f_+\|_p=\|g\|_{p'}^{1/(p-1)}.
\]
Since $p'-1/(p-1)=1$, division gives
\begin{equation}\label{eq:holder-first}
 \frac{\|Af_+\|_q}{\|f_+\|_p}
 \ge\frac{\|g\|_{p'}}{\|Af\|_q^{q-1}}.
\end{equation}

\textbf{2. Compare that pairing with the old ratio.}
Another use of adjointness followed by H\"older, now on $X$, gives
\[
 \|Af\|_q^q=\int_Y(Af)^q\,d\nu
          =\int_X fg\,d\mu
          \le\|f\|_p\|g\|_{p'}.
\]
Dividing by $\|f\|_p\|Af\|_q^{q-1}$ yields
\begin{equation}\label{eq:holder-second}
 \frac{\|g\|_{p'}}{\|Af\|_q^{q-1}}
 \ge\frac{\|Af\|_q}{\|f\|_p}.
\end{equation}
Combining the two comparisons proves monotonicity.

\textbf{3. Equality cannot hide in an intermediate inequality.}
If the first and last ratios agree, the middle quantity lies between
two equal numbers, so both comparisons are equalities. In particular,
H\"older on $X$ is equality. Because $f,g>0$ and $1<p<\infty$, its
precise equality condition is
\[
 f^p=a g^{p'}\quad\text{almost everywhere},\qquad a>0.
\]
Taking positive $p$th roots and using $p'/p=1/(p-1)$ gives
$f=a^{1/p}g^{1/(p-1)}=a^{1/p}f_+$. Thus the new function is a positive
multiple of the old function. Conversely, if $f_+=bf$ with $b>0$,
linearity of $A$ and homogeneity of both norms make the two ratios equal.
This proves the necessary and sufficient equality statement.
\end{proof}

\begin{proof}[Proof of Theorem~\ref{intro:lyapunov}]
Put $A=\mathcal A_m$, $f=\rho_K^m$, $p=n/m$ and $q=m+1$.
Since $K$ is a star body, $f$, $Af$ and
$g=A^*[(Af)^m]$ are continuous, bounded and bounded away from zero.
All norms and pairings in Lemma~\ref{lem:nonlinear-power} are
therefore finite and positive. There is no need here to assume a
bound for $A$ from the entire space $L^p$ into $L^q$.

By \eqref{eq:operator-normalized},
\[
 \rho_{\mathbb I_mK}^{m}
 =C_{n,m}^{m/(n-m)}g^{m/(n-m)}
 =C_{n,m}^{m/(n-m)}g^{1/(p-1)}.
\]
The positive scalar cancels in the ratio \eqref{intro:energy}.
Lemma~\ref{lem:nonlinear-power} gives
$\mathcal E_m(\mathbb I_mK)\ge\mathcal E_m(K)$, with equality
exactly when $\rho_{\mathbb I_mK}^m$ is proportional to $\rho_K^m$
almost everywhere. Both functions are continuous, so equality holds
everywhere; their positive $m$th roots give
$\mathbb I_mK=c_1K$. The converse follows from scale invariance.

If $K_j=\mathbb I_m^{\,j}K$ and $K_r=cK$, every $K_j$ remains a
star body by Lemma~\ref{lem:averaging}, and
\[
 \mathcal E_m(K_0)\le\mathcal E_m(K_1)\le\cdots\le
 \mathcal E_m(K_r)=\mathcal E_m(cK_0)=\mathcal E_m(K_0).
\]
All inequalities are equalities. The first one already implies that
$K$ is a fixed point up to dilation.
\end{proof}

\begin{proof}[Proof of Theorem~\ref{intro:main}]
Suppose $n\ge3$, $1\le m\le n-2$ and $\mathbb I_m^{\,r}K=cK$.
Theorem~\ref{intro:lyapunov} gives $\mathbb I_mK=c_1K$ for some
$c_1>0$. Theorem~\ref{thm:global-fixed} then forces $K$ to be an
origin-centered ball. Conversely, the image of every origin-centered
ball is another origin-centered ball, as computed in the proof of
Theorem~\ref{thm:global-fixed}. All its iterates are therefore positive
dilations of the original body.
\end{proof}

\begin{remark}[Scale invariance and the size of the iterates]
The geometric form of the Lyapunov functional is
\begin{equation}\label{eq:Lyapunov-geometry}
 \mathcal E_m(K)
 =\frac{(m+1)^{1/(m+1)}}{\omega_m}
   \left(\frac{\sigma_n}{n}\right)^{m/n}
   \frac{\Ph(K)^{1/(m+1)}}{V(K)^{m/n}}.
\end{equation}
If $\mathbb I_mK=c_1K$ and $\alpha=m^2/(n-m)$, homogeneity yields
\[
 \mathbb I_m^{\,r}K=c_1^{1+\alpha+\cdots+\alpha^{r-1}}K.
\]
These factors need not equal one. Theorem~\ref{intro:main} excludes
periodic changes of shape; it does not normalize the size of an orbit.
\end{remark}

\section{Local convergence near the ball}\label{sec:local}

In contrast with the global periodic-point classification, this section
is perturbative. We normalize by spherical mean, removing the constant
scaling direction, and use the spectrum computed in
Section~\ref{sec:regularity}.

\begin{proposition}[Local contraction in the stated Sobolev norm]
\label{prop:local-gap}
Let $2\le m\le n-2$ and $s>(n-1)/2$. With sphere probability $\sigma$,
set $H^s_0=\{h\in H^s:\int h\,\dd\sigma=0\}$ and
\[
 \mathcal T_m(\rho)
 =\frac{\mathcal P_m(\rho)}{\int\mathcal P_m(\rho)\,\dd\sigma}.
\]
On a neighborhood of $1$ in $H^s$, this map is smooth, with values
in $H^{s+m-1}$. Its derivative on $H^s_0$ is
\begin{equation}\label{eq:local-derivative}
 D\mathcal T_m(1)h=\frac{m^2}{n-m}N_mh,
 \qquad
 \|D\mathcal T_m(1)\|_{H^s_0\to H^s_0}=\frac{m}{n-1}<1.
\end{equation}
For each $\theta\in(m/(n-1),1)$ there is $\varepsilon>0$ such that
all $\rho$ with $\int\rho\,d\sigma=1$ and
$\|\rho-1\|_{H^s}<\varepsilon$ satisfy
\begin{equation}\label{eq:local-rate}
 \|\mathcal T_m^{\,j}(\rho)-1\|_{H^s}
 \le\theta^j\|\rho-1\|_{H^s},\qquad j\ge0.
\end{equation}
Here the exact operator norm refers to the spectral Sobolev norm
specified in Section~\ref{sec:regularity}.
\end{proposition}
\begin{proof}
\textit{Smoothness of the normalized map.}
Put $a=(m-1)/2$ and
$\mathcal B=H^{s+a}(\GG(n,m))\cap L^\infty(\GG(n,m))$ with its
sum norm. By the sphere embedding and the two separate estimates,
\[
 \|\mathcal A_m f\|_{\mathcal B}
 \le C\|f\|_{H^s}+\|f\|_\infty
 \le C'\|f\|_{H^s}.
\]
Thus $\mathcal A_m:H^s(\Sph^{n-1})\to\mathcal B$ is bounded linear.
Both $H^s(\Sph^{n-1})$ and $\mathcal B$ are Banach algebras,
the latter by Lemma~\ref{lem:powers}. Their integer-power maps are
continuous polynomials. Consequently
\[
 W(\rho):=\mathcal A_m^*[(\mathcal A_m(\rho^m))^m]
 :H^s(\Sph^{n-1})\longrightarrow H^{s+m-1}(\Sph^{n-1})
\]
is smooth and $W(1)=1$. Since $s+m-1>(n-1)/2$, its target is also a
Banach algebra. On a sufficiently small neighborhood of $1$,
the convergent binomial series for $(1+z)^{1/(n-m)}$ in that algebra
shows directly that $W(\rho)^{1/(n-m)}$ is smooth. The positive mean
is a smooth nonzero scalar there, so division by it preserves smoothness.
This proves smoothness with the stated gain of regularity.

\textit{Linearization and the spectral norm.}
The two integer powers and the final root give
\[
 DW(1)[h]=m^2N_mh,\qquad
 D\mathcal P_m(1)[h]=\alpha N_mh,\qquad
 \alpha=\frac{m^2}{n-m}.
\]
As $N_m1=1$ and $N_m$ is self-adjoint,
$\int N_mh\,\dd\sigma=\int h\,\dd\sigma$. Hence on the full space
\[
 D\mathcal T_m(1)[h]
 =\alpha\left(N_mh-\Bigl(\int h\,\dd\sigma\Bigr)1\right).
\]
On $H^s_0$ the constant term is zero. From \eqref{eq:eig},
\[
 \frac{\lambda_{2(k+1)}}{\lambda_{2k}}
 =\frac{(k+\frac12)(k+\frac{n-m}{2})}
 {(k+\frac m2)(k+\frac{n-1}{2})}<1,
\]
because denominator minus numerator equals
$(m-1)(k+n/4)>0$. Odd degrees have eigenvalue zero. Thus the
largest mean-zero eigenvalue is
\[
 \lambda_2=\frac{n-m}{m(n-1)},\qquad
 \alpha\lambda_2=\frac{m}{n-1}.
\]
Orthogonality for the spectral $H^s$ norm proves
\eqref{eq:local-derivative}.

\textit{The invariant neighborhood.}
Fix $\theta$ in the stated interval. Continuity of the derivative
allows a sufficiently small convex ball in $1+H^s_0$ on which all
restricted derivative norms are at most $\theta$. Shrink it so that
all radial functions are positive. For every normalized $\rho$ in this ball,
$\mathcal T_m(\rho)$ is normalized by definition. Thus derivatives along
mean-zero directions also take values in $H^s_0$.
The line segment $1+\tau(\rho-1)$ stays in the same convex ball, and
\[
 \mathcal T_m(\rho)-1
 =\int_0^1D\mathcal T_m(1+\tau(\rho-1))[\rho-1]\,d\tau.
\]
Taking the norm and using the uniform restricted derivative bound gives
\[
 \|\mathcal T_m(\rho)-1\|_{H^s}
 \le\theta\|\rho-1\|_{H^s}.
\]
The ball is therefore invariant, and iteration proves
\eqref{eq:local-rate}. The sharp linear factor $m/(n-1)$ is not
asserted to be an exact nonlinear contraction factor on a fixed
neighborhood.
\end{proof}

\begin{remark}[The degree-two multiplier]
The first nonconstant multiplier can also be computed without special
functions. If $h(v)=v^TSv$, where $S$ is symmetric and
$\operatorname{tr}S=0$, then
\[
 \mathcal A_mh(\xi)=\frac1m\operatorname{tr}(P_\xi S).
\]
For a uniformly distributed $m$-plane containing $u$,
\[
 \mathbb E P_\xi=uu^T+\frac{m-1}{n-1}(I-uu^T).
\]
Taking the trace against $S$ gives
$N_mh=(n-m)h/(m(n-1))$, and hence the nonlinear multiplier
$m/(n-1)$. For $n=4,m=2$, the full formula reduces to
$\lambda_{2k}=1/(2k+1)$ and this multiplier is $2/3$.
\end{remark}

\begin{proposition}[Second variation of the Lyapunov ratio at the ball]
\label{prop:energy-hessian}
Let $2\le m\le n-2$, $p=n/m$ and $q=m+1$. Set
$E(f)=\|\mathcal A_m f\|_q/\|f\|_p$ using probability measures.
For every bounded real $h$ with $\int h\,d\sigma=0$,
\begin{equation}\label{eq:energy-hessian}
 \left.\frac{d^2}{d\epsilon^2}\right|_{0} E(1+\epsilon h)
 =m\langle h,N_mh\rangle-\frac{n-m}{m}\|h\|_2^2.
\end{equation}
Consequently this is at most
\begin{equation}\label{eq:energy-hessian-gap}
 -\frac{(n-m)(n-m-1)}{m(n-1)}\|h\|_2^2,
\end{equation}
which is strictly negative for $h\ne0$.
\end{proposition}
\begin{proof}
For $|\epsilon|\|h\|_\infty<1/2$, both $1+\epsilon h$ and
$1+\epsilon\mathcal A_mh$ lie in $[1/2,3/2]$. Taylor's formula
for the real powers $p,q$, followed by that for their positive roots,
is uniform on this interval. Since both first moments vanish,
\begin{align*}
 \|1+\epsilon h\|_p
 &=1+\frac{p-1}{2}\epsilon^2\|h\|_2^2+O_h(\epsilon^3),\\
 \|1+\epsilon\mathcal A_mh\|_q
 &=1+\frac{q-1}{2}\epsilon^2\|\mathcal A_mh\|_2^2
       +O_h(\epsilon^3).
\end{align*}
Dividing the two expansions and using
$\|\mathcal A_mh\|_2^2=\langle h,N_mh\rangle$ proves
\eqref{eq:energy-hessian}. The largest eigenvalue of $N_m$ on the
mean-zero subspace is $(n-m)/(m(n-1))$, while odd harmonics are
annihilated. Substitution gives \eqref{eq:energy-hessian-gap}.
\end{proof}

The increasing Lyapunov functional and local attraction are therefore
consistent: nonconstant perturbations lower its value to second order,
and the local iteration raises it toward the value at the ball.
Neither the second variation nor the local contraction replaces the
global argument of Section~\ref{sec:weighted-strict}. No uniform
compactness or convergence assertion for arbitrary normalized orbits
is made here.

\appendix
\section{A singular-integral representation}\label{app:singular}

This representation is independent of the graph formula used in the
main proof. Let $\Delta_j$ denote the $j$-dimensional parallelotope
volume. The linear Blaschke--Petkantschin formula, with Haar
probability on $G(n,m)$, reads
\begin{equation}\label{eq:BP}
 \int_{(\R^n)^m}g(\mathbf x)\,d\mathbf x
 =b_{n,m}\int_{G(n,m)}\int_{\xi^m}
    \Delta_m(\mathbf x)^{n-m}g(\mathbf x)\,d\mathbf x\,d\nu_{n,m}(\xi),
\end{equation}
where
\[
 b_{n,m}=\frac{\sigma_n\sigma_{n-1}\cdots\sigma_{n-m+1}}
                  {\sigma_m\sigma_{m-1}\cdots\sigma_1}.
\]
Here the determinant is a parallelotope volume, not a simplex volume;
this convention fixes the constant. See \cite{Gardner2006} for the
integral-geometric framework.

\begin{proposition}[Blaschke--Petkantschin concentration formula]\label{prop:first}
If $K\subset\R^n$ is compact and $\HH^m(\partial K\cap\xi)=0$
for $\nu_{n,m}$-almost every $\xi$, then
\begin{equation}\label{eq:first}
 \Ph(K)=\lim_{p\downarrow m-n}
 \frac{p+n-m}{(m+1)b_{n,m}\sigma_{n-m}}
 \int_{K^{m+1}}\Delta_{m+1}(x_1,\ldots,x_{m+1})^p
 \dd x_1\cdots\dd x_{m+1}.
\end{equation}
\end{proposition}
\begin{proof}
Take $R\ge1$ with $K\subset B_2^n(R)$, put $k=n-m$, and write
$p=-k+\varepsilon$, $0<\varepsilon<1$. On full-rank tuples in $\xi^m$,
\[
 \Delta_{m+1}(x_1,\ldots,x_m,z)
 =\Delta_m(x_1,\ldots,x_m)|P_{\xi^\perp}z|.
\]
Applying \eqref{eq:BP} to the first $m$ variables, as an identity
for nonnegative measurable functions, transforms the integral into
\[
 b_{n,m}\int_{\GG(n,m)}\int_{(K\cap\xi)^m}
 \Delta_m(\mathbf x)^\varepsilon
 \int_K|P_{\xi^\perp}z|^{-k+\varepsilon}\,\dd z
 \,\dd\mathbf x\,\dd\nu_{n,m}(\xi).
\]
Let
\[
 T_\varepsilon(\xi)=\frac{\varepsilon}{\sigma_k}
       \int_K|P_{\xi^\perp}z|^{-k+\varepsilon}\,\dd z.
\]
Split $z=y+w\in\xi\oplus\xi^\perp$ and use polar coordinates in
$w$. For $y\notin\partial K$, the function $w\mapsto\one_K(y+w)$
is constant near $w=0$. Its integral against
$\varepsilon r^{\varepsilon-1}\dd r$ on $[0,R]$, averaged over the
unit sphere in $\xi^\perp$, therefore tends to $\one_K(y)$:
the mass of $[0,R]$ is $R^\varepsilon\to1$, whereas that of
$[\eta,R]$ is $R^\varepsilon-\eta^\varepsilon\to0$ for every $\eta>0$.
The assumed boundary condition and dominated convergence in
$y\in\xi\cap B_2^n(R)$ imply
$T_\varepsilon(\xi)\to\vol_m(K\cap\xi)$ for almost every $\xi$.
The domination is uniform, since
\[
 0\le T_\varepsilon(\xi)\le\omega_mR^{m+\varepsilon}
 \le\omega_mR^{m+1}.
\]
On $(\xi\cap B_2^n(R))^m$ we also have
$\Delta_m(\mathbf x)^\varepsilon\le R^m$. Thus the normalized
integrand is bounded by $\omega_mR^{2m+1}$ on a bundle of total
measure $(\omega_mR^m)^m$. Rank-deficient tuples have fiber Lebesgue
measure zero, and $\Delta_m^\varepsilon\to1$ on full-rank tuples.
A final dominated-convergence application gives
$(m+1)^{-1}\int\vol_m(K\cap\xi)^{m+1}\dd\nu$, as claimed.
\end{proof}

For a star body, the boundary hypothesis is automatic: on every
$m$-plane, polar coordinates give zero $m$-dimensional measure to the
radial graph $\{\rho_K(u)u:u\in\Sph^{n-1}\cap\xi\}$.
The graph representation \eqref{eq:weighted-graph}, however, avoids this
limit altogether and is the representation used for symmetrization.

\section{An elementary rigidity argument in a restricted range}\label{app:maxmin}

The following proof uses neither regularity nor symmetrization. It is
included to exhibit directly the role of the homogeneity threshold
$n=m(m+1)$; it is not needed in the proof of the main theorem.

\begin{proposition}[Global rigidity in the subcritical and critical regimes]
\label{prop:maxmin}
Let $K$ be an arbitrary star body in $\R^n$, let $1\le m\le n-2$,
and suppose $n\ge m(m+1)$. If $\mathbb I_mK=cK$ for some $c>0$,
then $K$ is a Euclidean ball centered at the origin.
\end{proposition}
\begin{proof}
Set
\[
 r=\min_{\Sph^{n-1}}\rho_K,
 \qquad R=\max_{\Sph^{n-1}}\rho_K,
 \qquad A=\frac{\sigma_m}{\sigma_n}\omega_m^m.
\]
Since $rB_2^n\subset K\subset RB_2^n$,
$\omega_m r^m\le\vol_m(K\cap\xi)\le\omega_m R^m$
for every $m$-plane $\xi$. The positivity of the dual Radon
transform gives
\[
 Ar^{m^2}\le\rho_{\mathbb I_mK}(u)^{n-m}
 \le AR^{m^2}\qquad\text{for all }u\in\Sph^{n-1}.
\]
Evaluate the lower bound at a minimizer of $\rho_K$ and the upper
bound at a maximizer. Under the fixed-point equation,
\[
 c^{n-m}\ge A r^{m^2-(n-m)},\qquad
 c^{n-m}\le A R^{m^2-(n-m)}.
\]
If $n>m(m+1)$, the common exponent is negative and hence these
two inequalities imply $r\ge R$, forcing $r=R$.

It remains to treat $n=m(m+1)$. As $m\le n-2$ and $n\ge3$, here
$m\ge2$. If $r<R$, continuity of $\rho_K$ yields a nonempty open
spherical cap on which $\rho_K>r$. Fix any direction $u$. Choose a point $v$ in that cap not parallel to
$u$; this is possible because a nonempty spherical open set is not
contained in two points. Since $m\ge2$, extend $\operatorname{span}\{u,v\}$
to an $m$-plane $\xi_0$. The cap meets the unit sphere in $\xi_0$ in a
nonempty open subset, which has positive $(m-1)$-dimensional measure.
Thus $|K\cap\xi_0|_m>\omega_mr^m$. Section volume depends continuously
on $\xi$ (rotate the sphere of integration and use uniform continuity
of $\rho^m$). This strict inequality holds on a nonempty open subset
of the incidence Grassmannian of planes through $u$. Its invariant
probability has full support, so that subset has positive measure.
The remaining sections satisfy the weak lower bound. Consequently,
\[
 \mathcal R_m^*\bigl[\vol_m(K\cap\xi)^m\bigr](u)>Ar^{m^2}.
\]
Likewise, a cap on which $\rho_K<R$ gives the strict upper bound
$\mathcal R_m^*[\vol_m(K\cap\xi)^m](u)<AR^{m^2}$.
At a radial minimizer and maximizer respectively these inequalities
imply both $c^{n-m}>A$ and $c^{n-m}<A$, a contradiction.
Thus $r=R$ also at the critical exponent.
\end{proof}

\begingroup
\small
\setlength{\bibsep}{4pt}

\endgroup

\end{document}